\documentclass[11pt,leqno]{article}
\newtheorem{thm}{Theorem}[section]
\newtheorem{lem}[thm]{Lemma}

\newtheorem{prop}[thm]{Proposition}

\newtheorem{definition}{Definition}

\newtheorem{rem}{Remark}

\renewcommand{\P}{{\bf P}}
\newcommand{\R}{{\bf R}}
\newcommand{\C}{{\bf C}}
\newcommand{\Z}{{\bf Z}}

\newcommand{\ELL}{{\bf  L}}

\newcommand{\bF}{{\bar{F}}}

\newcommand{\cI}{{\cal I}}
\newcommand{\cG}{{\cal G}}
\newcommand{\cL}{{\cal L}}
\newcommand{\cO}{{\cal O}}
\newcommand{\cR}{{\cal R}}

\newcommand{\al}{{\alpha}}
\newcommand{\be}{{\beta}}
\newcommand{\ep}{{\varepsilon}}
\newcommand{\lm}{{\lambda}}
\newcommand{\vf}{{\varphi}}
\newcommand{\gm}{{\gamma}}

\newcommand{\Aut}{{\mathop{\mathrm{Aut}}\nolimits}}
\newcommand{\GL}{{\mathop{\mathrm{GL}}\nolimits}}
\newcommand{\SL}{{\mathop{\mathrm{SL}}\nolimits}}
\newcommand{\PGL}{{\mathop{\mathrm{PGL}}\nolimits}}
\newcommand{\PSL}{{\mathop{\mathrm{PSL}}\nolimits}}
\newcommand{\Uni}{{\mathop{\mathrm{U}}\nolimits}}
\newcommand{\Grass}{\mathop{\mathrm{Gr(2,4)}}\nolimits}
\newcommand{\Grassm}{\mathop{\mathrm{Gr(m,2m)}}\nolimits}
\newcommand{\Grassn}{\mathop{\mathrm{Gr(n+1,2n+2)}}\nolimits}
\newcommand{\Gr}{\mathop{\mathrm{Gr}}\nolimits}
\newcommand{\OG}{{\Omega/\Gamma}}

\newcommand{\wg}{\wedge}
\newcommand{\Wg}{\bigwedge}
\newcommand{\rank}{\mathop{\mathrm{rank}}\nolimits}
\newcommand{\Ker}{\mathop{\mathrm{Ker}}\nolimits}
\renewcommand{\Im}{\mathop{\mathrm{Im}}\nolimits}
\newcommand{\Int}{\mathop{\mathrm{Int}}\nolimits}
\newcommand{\pd}{{\partial}}
\newcommand{\Lam}{\mathop{\mathrm{\Lambda^m(\C^{2m})}}\nolimits}

\newcommand{\bDelta}{{\bar{\Delta}}}

\newcommand{\ssig}{{\{\sigma_\nu\}}}
\newcommand{\tssig}{{\{\tilde{\sigma}_\nu\}}}
\newcommand{\tsigma}{{\tilde{\sigma}}}
\newcommand{\tx}{{\tilde{x}}}
\newcommand{\tg}{{\tilde{g}}}
\newcommand{\ttau}{{\tilde{\tau}}}
\newcommand{\tGamma}{{\tilde{\Gamma}}}

\newcommand{\hsigma}{{\hat{\sigma}}}
\newcommand{\hsig}{{\{\hat\sigma_\nu\}}}

\newcommand{\dist}{\mathop{\mathrm{distance}}\nolimits}

\newcommand{\Mt}[1]{{\mathrm{M}_{#1}(\bf C)}}
\newcommand{\MAT}[4]{\left(\begin{array}{cc} #1 & #2 \\ #3 & #4 \end{array}\right)}
\newcommand{\Frac}[2]{\displaystyle{\frac{#1}{#2}}}
\newcommand{\qed}{\hbox{\rule[0pt]{3pt}{6pt}}}

\title{On odd dimensional complex analytic Kleinian groups}
\author{Masahide Kato}
\date{}

\begin{document}

\maketitle

\begin{abstract}
We shall explain here an idea to generalize classical 
complex analytic Kleinian  group theory to any odd dimensional cases.
For a certain class of discrete subgroups of $\PGL_{2n+1}(\C)$ acting on 
$\P^{2n+1}$, we can define their domains of discontinuity 
in a canonical manner, 
regarding an $n$-dimensional projective linear subspace in $\P^{2n+1}$
as a point, like a point in the classical $1$-dimensional case. 
Many interesting (compact) non-K\"ahler manifolds appear systematically 
as the canonical quotients of the domains. 
In the last section, we shall give some examples.
\footnote{keywords : compact non-K\"ahler manifold, 
projective structure, Kleinian group theory}

\end{abstract}

\section{Introduction}\label{s:1}
The theory of discrete subgroups of $\PGL_2(\C)$ 
has a long history.
Let $\Gamma$ be a discrete subgroup of $\PGL_2(\C)$.
We say that the action of $\Gamma$ 
at a point $z \in \P^1$ is {\em discontinuous}, 
if there is a neighborhood $W$ of $z$ 
such that $\gamma(W)\cap W =\emptyset$ for 
all but finitely many $\gamma\in \Gamma$.
Following B. Maskit \cite{Maskit}, we call  a 
subgroup $\Gamma \subset \PGL_2(\C)$ 
whose action is discontinuous at some point $z \in \P^1$ by  
{\em Kleinian group}. 

Let $\Gamma\subset \PGL_2(\C)$ be a Kleinian group. 
The set $\Omega(\Gamma)$ of points $z \in \P^1$  at which $\Gamma$ acts 
discontinuously is called the {\em set of discontinuity} of $\Gamma$. 
The set $\Omega(\Gamma)$ is a $\Gamma$-invariant open subset in $\P^1$ on 
which $\Gamma$ acts properly discontinuously. 
The geometry of the quotient space $\Omega(\Gamma)/\Gamma$ is one of 
the main themes in the classical Kleinian group theory. 

If we seek for a higher dimensional version of the Kleinian group theory,  
we must first define the set of discontinuity for a given discrete subgroup. 
Let $n\geq 2$. 
Take a discrete subgroup $\Gamma \subset \PGL_{n+1}(\C)$ acting of $\P^n$.
Consider, as above,  the set $\Omega(\Gamma)$ of points $z \in \P^n$  
at which $\Gamma$ acts discontinuously.
Then it is true that $\Gamma$ acts on $\Omega(\Gamma)$, 
but the action is not properly discontinuous 
in general. 
Therefore, we must find another definition of the set of discontinuity to get a good
quotient space.

In this paper, we consider 
a class of discrete subgroups in $\PGL_{2n+2}(\C)$ ($n\geq 1$), i.e, 
a class of {\it type $\ELL$ groups} (Definition \ref{DefOfL}). 
A type $\ELL$ group $\Gamma$ has the non-empty set of discontinuity 
$\Omega(\Gamma)$, which is defined in a canonical manner. 
The set $\Omega(\Gamma)$ contains  a subdomain 
$W \subset \P^{2n+1}$, which is biholomorphic to    
\begin{equation}\label{U}
\{ z \in \P^{2n+1}
 : |z_0|^2 +\cdots + |z_n|^2 < |z_{n+1}|^2 +\cdots + |z_{2n+1}|^2 \} 
\end{equation}
and satisfies
\begin{equation} 
\gamma(W) \cap W = \emptyset \mbox{ for any } \gamma \in \Gamma\setminus\{1\},
\end{equation}
where
$z= [z_0:\dots:z_n:z_{n+1}:\dots:z_{2n+1}]$.

\bigskip
For type $\ELL$ groups,  
$n$-dimensional projective linear subspaces in $\P^{2n+1}$ play 
the same role as points do in one dimensional Kleinian group theory. 
In the following, an $n$-dimensional projective linear subspace 
is called an {\em $n$-plane} for short.
The paper is organized as follows. 

In section $1$, we shall make some preparations on the Grassmannian
$\Grassm$ of $m$-dimensional subspaces in $\C^{2m}$. 
As is well-known, $\Grassm$ can be embedded into 
the projective space $\P^N$, $N={{_{2m}C}_m -1}$, 
by Pl\"ucker coordinates. 
We remark that the embedded $\Grassm$ is 
contained in a $\PGL_m(\C)$-invariant hyperquadric (Proposition \ref{GinQ}). 
This fact plays an important role in studying limit sets of type $\ELL$ groups.
In section 2, we study some convergence properties of infinite sequences of 
projective transformations.
In section 3, we define the set of discontinuity $\Omega(\Gamma)$ 
for a type $\ELL$ group $\Gamma$ (Definition \ref{ProperlyDiscon}, 
Definition \ref{DefOfL}),
and show that the action of $\Gamma$ on $\Omega(\Gamma)$ is 
properly discontinuous (Theorem \ref{CanQ}).  
Hence the quotient $\Omega(\Gamma)/\Gamma$ becomes a good space.
A domain $\Omega \subset \P^{2n+1}$ is said to 
be {\em large}\footnote{The definition of the term "large"
is different from \cite{Lar}, where a domain $\Omega\subset\P^{2n+1}$ 
is said to be large if the $4n$-dimensional Hausdorff measure of its complement vanishes},
if $\Omega$ contains an $n$-plane. 
Any holomorphic automorphism of a large domain 
extends to an element of $\PGL_{2n+2}(\C)$ (Ivashkovich \cite{Iv}). 
Using this, we show in section 4, 
that a large domain which covers a compact manifold 
is a connected component of $\Omega(\Gamma)$ of some type $\ELL$ 
group $\Gamma$  (Theorem \ref{component}).
This may justify our definition of $\Omega(\Gamma)$.
There are many groups of type $\ELL$. As an example, 
we explain briefly an analogue of Klein combinations and 
handle attachments in section 5.
See also \cite{line} on this topic.
In section 6, an analogue of the Ford region is 
defined. We prove that this region 
gives a fundamental set of a type $\ELL$ group under some additional conditions. 
In section 7, we shall give examples of type $\ELL$ groups and their quotient 
spaces $\Omega(\Gamma)/\Gamma$. 

\tableofcontents
\bigskip
{\bf Notation}
\begin{itemize}
\item $M(p\times q,\C)$: the set of matrices of size $p\times q$ 
with coefficients in $\C$.
\item $M_p(\C)$: the set of matrices of size $p\times p$ 
with coefficients in $\C$.
\end{itemize}

\section{The Grassmannian $\Grassm$}
Let $\Grassm$, $m\geq 2$, be the Grassmannian 
of the $m$-dimensional subspaces in $\C^{2m}$.
The aim of the section is to show that $\Grassm$ is embedded in 
a quadric hypersurface in a big projective space by Pl\"ucker coordinates. 
This is a well-known fact. But since this is important for the later argument, 
we will explain it here for the reader's convenience.

Let $\{e_1,\dots, e_{2m}\}$ be a basis of $\C^{2m}$, and 
$\cI$ the set of multiindices 
\[ I = \{i_1,\dots, i_m\} \subset \{1,\dots, 2m\}, 
\ \ \  i_1 < \dots < i_m,
\]
of cardinality $m$. 
In the set of multiindices $\cI=\{I\}$, 
we introduce the lexicographic order.
Namely, for multiindices   
$I=\{i_1,\dots, i_m\}$, 
$J=\{j_1,\dots, j_m\}$, $I \neq J$,  
we write $I < J$, if $i_\mu < j_\mu$ for 
$\mu = \min\{\lm : i_\lm \neq j_\lm\}$.
We put
\begin{equation}\nonumber
\delta_{JK}  = \delta^{KJ}=
\left\{\begin{array}{cl}
(-1)^\nu, & J \cap K = \emptyset, \\
0,             & J \cap K \neq \emptyset,
\end{array}\right. 
\end{equation}
where $\nu = \# \{ (p, q) \in J \times K : p > q \}$, and
\begin{equation}\nonumber
\delta_J^I  =  
\left\{\begin{array}{cl}
1,     & I = J, \\
0,     & I \neq J.
\end{array}\right.
\end{equation}
Then we have $\delta_{IJ}=(-1)^m\delta_{JI}$, $\delta^{IJ}=(-1)^m\delta^{JI}$, and 
\begin{equation}\nonumber
\delta_{IJ}\delta^{JK} = \left\{\begin{array}{cc} 1, & I=K, \\ 0, & I\neq K. \end{array}\right. \ \ \mbox{(Einstein's convention)}
\end{equation}
As a basis of $\Lambda^m(\C^{2m})$ of the space of $m$-vectors, 
we use  $\{ e_I \}_I$, where
\[ e_I = e_{i_1}\wg \dots \wg e_{i_m}, \ \ \ I=\{i_1,\dots, i_m\} \in \cI. \]
Then we have
\begin{equation} \label{JK}
e_J\wg e_K = \delta_{JK}e_1\wg \dots \wg e_{2m}. 
\end{equation}
Any $w \in \Lam$ is written uniquely as a linear combination over $\C$,
\[ w = w^I e_I, \ \ \ w^I \in \C, \ \ \ \mbox{(Einstein's convention)}. \]
If $w\neq 0$, it determines a point $[w^I] \in \P^{N}$,
regarding $\{w^I\}_I$ as a homogeneous coordinates, where 
$N={_{2m}C}_m -1$. 

Let $X$ be an $m$-dimensional subspace in $\C^{2m}$ spanned by 
$2m$-vectors $\{x_1,\dots, x_m\}$. Then $X$ corresponds
to the $m$-vector $\hat{X}= x_1\wg \dots \wg x_m \in \Grassm$. 
Letting $x_j = x_j^k e_k$, we have 
\begin{equation}\label{plucker}
\hat{X} 
 =  x_1^{k_1}e_{k_1} \wg \dots \wg x_m^{k_m}e_{k_m}
 =  X^Ke_K,
\end{equation}
where 
\[ 
X^K =\det \left(
\begin{array}{ccc}
x_1^{k_1} & \cdots & x_m^{k_1} \\
\vdots              &            &  \vdots \\
x_1^{k_m} & \cdots & x_m^{k_m} 
\end{array}
\right).
\]
The set of numbers $\{X^K\}_{K \in \cI}$ determines the point 
$[x^K] \in \P^N$, 
which is the Pl\"ucker coordinates of the vector subspace $X$. 

\bigskip
Let $A \in M_{2m}(\C)$ be any element.  Put $Ae_j = a^k_je_k$.
Then we have  
\[
 Ae_J  =  Ae_{j_1}\wg \dots \wg Ae_{j_m} 
          =  a^{k_1}_{j_1}e_{k_1}\wg \dots \wg a^{k_m}_{j_m}e_{k_m}
          =  A_J^K e_K,
\]
where 
\[
A_J^K =\det \left(
\begin{array}{ccc}
a_{j_1}^{k_1} & \cdots & a_{j_m}^{k_1} \\
\vdots              &            &  \vdots \\
a_{j_1}^{k_m} & \cdots & a_{j_m}^{k_m} 
\end{array}
\right).
\]
Hence, for $J, K \in \cI$ with $J\cap K=\emptyset$, we have
\[
A(e_J\wg e_K)   
         =  A_J^L e_L\wg A_K^M e_M 
          = \delta_{LM}A^L_J A^M_K e_1\wg\dots\wg e_{2m}.
\]
On the other hand, we have
\[
A(e_J\wg e_K) =  \delta_{JK}\det(A)e_1\wg\dots\wg e_{2m}
\]
Thus we have
\begin{equation}\label{AA}
\delta_{LM}A^L_J A^M_K=\delta_{JK}\det A.
\end{equation}
Define a bilinear form $Q(z,w)$ on $\C^{N+1}$ by
\begin{equation}\nonumber 
Q(z, w) = \delta_{JK}z^Jw^K, \ \ \ z=(z^J), \ \ w=(w^K). 
\end{equation}
Put 
\[ \hat A z = (A^K_Iz^I), \ \ \ z=(z^I). \]
Then, by (\ref{AA}), we have 
that
\begin{equation}\nonumber 
Q(\hat{A}z, \hat{A}w) = (\det A)Q(z,w).
\end{equation}
For $\hat{A}$, we define $\hat{A}^* \in M_{N+1}(\C)$ by
\begin{equation}\label{Astar} 
(\hat{A}^*)^I_J = \delta^{IK}\delta_{LJ}A^{L}_{K}
\end{equation}
Then we have 
\begin{eqnarray}  
\label{invariant} Q(\hat{A}z, w) & = & Q(z, \hat{A}^*w), \\ 
\label{inverse} \hat{A}^*\hat{A} & = & (\det A) I_{N+1}.
\end{eqnarray}
\begin{prop}\label{QXY} 
Let $X, Y \subset \C^{2m}$ be 
$m$-dimensional vector subspaces.
Put $X=X^Ke_K$ and $Y=Y^Ke_K$.
Then $\dim (X\cap Y)\geq 1$ holds if and only if   
\[ Q\left((X^K), (Y^K)\right)=0. \]
In particular, the equation 
\[ Q\left((X^K), (X^K)\right)=0 \]
holds for any $m$-dimensional subspace $X$ of $\C^{2m}$. 
\end{prop}
{\bf Proof} \  Clear by (\ref{plucker}), (\ref{JK}). \ \qed

\bigskip
We apply the above argument to the case $\Grassn$. 
Set $N=_{2n+2}C_{n+1}-1$. Then, by Proposition \ref{QXY}, we have easily
\begin{prop}\label{PQXY} 
Let $Q(z,w)$ be the quadratic form defined by 
\[ Q(z,w) = \delta_{JK}z^Jw^K \]
defined on $\C^{N+1}\times \C^{N+1}$.
Let $\ell_1, \ell_2$ be $n$-planes in $\P^{2n+1}$,
and $[x^I]$, $[y^I]$ be their Pl\"ucker coordinates.
Then $\ell_1$ and $\ell_2$ intersects
if and only if  the equality 
\[ Q([x^I],[y^I])=0 \]
holds. 
In particular, for an $n$-plane with the Pl\"ucker coordinates
$[x^I]$,  we have 
\begin{equation}\label{Qxx}
 Q([x^I],[x^I])=0. 
\end{equation}
\end{prop}
We have also the following
\begin{prop}\label{GinQ} 
The quadric hypersurface $Q=\{Q(z,z)=0\}$ 
in $\P^N$, 
is invariant by the image group of the group representation
\[
\rho \ : \ \PGL_{2n+2}(\C) \to \PGL_{N+1}(\C), \ \ \ \rho(A)=\hat A,
\]
and the Grassmannian $\Grassn$ is contained in $Q$.
\end{prop}

\section{Limit of projective transformations}\label{foundation}  
Let $N\geq 1$ and $\Gamma$ a discrete infinite subgroup of $\PGL_{N+1}(\C)$
which acts on the projective space $\P^N$.
Consider an infinite sequence $\ssig$ of elements of $\Gamma$. 
Let $\tilde\sigma_\nu \ \in\GL_{N+1}(\C)$ be a representative 
of $\sigma_\nu$ such that $|\tilde\sigma_\nu|=1$, 
where, for a matrix $A= (a_{jk})$ of size $N+1$, 
we put $|A| = \max_{0\leq j,k\leq N}|a_{jk}|$. 
We say that $\ssig$ is a {\em normal sequence}, 
if the following conditions are satisfied.
\begin{enumerate}
\item The sequence $\ssig$ consists of distinct elements of $\Gamma$.
\item The sequence of matrices $\tssig$ can be chosen to be 
convergent to a matrix $\tilde\sigma \in M_{N+1}(\C)$. 
\end{enumerate}
The projective linear subspace defined by the image of  
the linear map $\tilde\sigma :\C^{N+1} \to \C^{N+1}$
is called the {\em limit image} of the normal sequence $\ssig$ 
and denoted by $I(\ssig)$.
Similarly the projective linear subspace defined by 
the kernel of $\tilde\sigma$ is called the {\em limit kernel} of 
$\ssig$ and denoted by $K(\ssig)$.
Here $r = \rank \tilde\sigma$ is called the rank of the normal sequence. 
Note that $I(\ssig)$, $K(\ssig)$, and $r$ are determined 
independently of the choice of representatives $\tilde\sigma_\nu$. 
Obviously, we have
$\dim I(\ssig)=r-1$ and $\dim K(\ssig)=N-r$.

\begin{thm}(\cite[Satz 2]{My})\label{Myr1}. 
Let $\{\sigma_\nu\} \subset \Gamma$ 
be a normal sequence. 
Suppose that the sequence of representatives $\{\tsigma_\nu\}$
converges to $\tsigma : \C^{N+1} \to \C^{N+1}$. 
Let $I$ be its limit image and $K$ the limit kernel.  
Then the sequence $\{\sigma_\nu\}$ converges uniformly on compacts in 
$\P^N\setminus K$ to 
the projection $\sigma : \P^N \setminus K \to I$ defined by $\tsigma$. 
\end{thm}

\begin{thm} \label{IinQ}
Let $\{\sigma_\nu\}_\nu \subset \Gamma$ be 
a normal sequence such that the sequence $\{\hsigma_\nu\}_\nu$ 
is also normal. Then the limit image of $\{\hsigma_\nu\}_\nu$
is contained in $Q$, and the limit kernel coincides with the orthogonal 
subspace (with respect to $Q(z,z)$) 
of the limit image of $\{{\hsigma_\nu}^{-1}\}_\nu$.
\end{thm}
{\bf Proof} \  
Let $S_\nu \in \GL_{N+1}(\C)$ be a  representative of $\hsigma_\nu$.
We can assume that $|S_\nu| = 1$, and that the sequence $\{S_\nu\}$ 
converges to $S \in M_{N+1}(\C)$.
Since $\Gamma$ is discrete, we have $\det S = \lim_\nu \det S_\nu=0$.
Therefore we have
\begin{equation}\label{QQ} 
Q(S_\nu z, S_\nu z) = (\det S_\nu) Q(z, z), \ \ \ \mbox{and}
\ \ \ Q(Sz, Sz) = 0, 
\end{equation}
by Proposition \ref{GinQ}. 
Hence the limit image of $\{\hsigma_\nu\}$ is contained in $Q$. Since 
\begin{eqnarray*}
(\Im S^*)^\perp 
& = & \{ z \in \C^{N+1} : Q(z, S^*w) = 0 \ \forall w \in \C^{N+1} \} \\
& = & \{ z \in \C^{N+1} : Q(Sz, w) = 0 \ \forall w \in \C^{N+1} \} \\
& = & \Ker S,
\end{eqnarray*}
we have 
\begin{equation}\label{KpI}
\Ker S =(\Im S^*)^\perp.
\end{equation} 
By (\ref{inverse}), we see that the projection $\P^N \dots\to \P^N$ 
defined by $S^*$ is the limit of
the normal sequence $\{\hsigma_\nu^{-1}\}_\nu$.
Thus we have the theorem. 
\ \qed

\section{Discontinuous groups in the projective $(2n+1)$-space}
From now on, {\em we assume that
$\Gamma\subset \PGL_{2n+2}(\C)$ is of type $\ELL$}, 
if not stated otherwise explicitly. This is our higher dimensional 
{\em complex analytic} analogue of Kleinian groups. 
Put $N={_{2n+2}}C_{n+1}-1$ and 
$\cG=\Gr(n+1,2n+2)$. 
We shall say, from now on, that 
a sequence $\ssig \subset \Gamma$ is {\em normal},
if not only the original sequence $\ssig$ is normal
but also is the corresponding sequence $\{\hat\sigma_\nu\}$, 
$\hat\sigma_\nu=\rho(\sigma_\nu)$, of $\PGL_{N+1}(\C)$. 
Thus a normal sequence $\ssig \subset \Gamma$ 
defines also $I(\hsig)$ and $K(\hsig)$ in $\P^N$.
Note that any normal sequence 
in the old sense contains a subsequence which is normal in 
the new one. 

\smallskip
\begin{definition}\label{Limitline} An $n$-plane $\ell$ in 
$\P^{2n+1}$ 
is called a {\em limit $n$-plane} of $\Gamma$, 
if there is a normal sequence $\ssig$ of $\Gamma$ with 
$\hat\ell \in \cG\cap I(\hsig)$.
\end{definition}

Let $\cL(\Gamma)\subset \cG$ denote 
the set of points which correspond to limit $n$-planes of $\Gamma$.

\begin{definition}\label{LimitSet} The union 
\[ \Lambda(\Gamma) = \bigcup_{\hat\ell\in \cL(\Gamma)}|\ell| \]
of the support of limit $n$-planes of $\Gamma$ is called 
{\em the limit set} of $\Gamma$.
\end{definition}  
Here we indicate by $|\ell|$ the support of an $n$-plane 
$\ell$ in $\P^{2n+1}$ in order to express explicitly the set of
points on the $n$-plane. 
\begin{definition}\label{ProperlyDiscon}
The set 
\[ \Omega(\Gamma) = \P^{2n+1}\setminus\Lambda(\Gamma) \]
is called the set of discontinuity of the group $\Gamma$. 
\end{definition}

\begin{definition}\label{DefLarge} A domain $\Omega$ 
in $\P^{2n+1}$ is said to be {\em large}, if $\Omega$ contains an $n$-plane.
\end{definition}

There are examples of $\Gamma$ with non-empty $\Omega(\Gamma)$, 
but which contains no $n$-planes. For example in the case $\P^3$, 
let $\Gamma$ be the infinite cyclic group generated by 
$\sigma=\MAT{I}{0}{A}{I} \in \PGL_4(\C)$, 
$A=\MAT{1}{0}{0}{0}$. 
Then $\Omega(\Gamma)=\P^3\setminus\{z_0=0\}$.
Thus we define type $\ELL$ groups as follows.

\bigskip
\begin{definition}\label{DefOfL} 
A discrete subgroup in $\PGL_{2n+2}(\C)$ is said to be of
type $\ELL$, if $\Omega(\Gamma)$ contains a large domain. 
\end{definition}
 
\begin{lem}
$\cG\cap I(\hsig)$ consists of a single point 
for any normal sequence $\ssig$ in $\Gamma$.
\end{lem}
{\bf Proof} \   Let $\{\sigma_\nu\}_\nu$ be any normal sequence 
in $\Gamma\subset\PGL_{2n+2}(\C)$, and $\{S_\nu\}_\nu \subset \GL_{N+1}(\C)$ be any 
convergent sequence of representatives of $\{\hsigma_\nu\}_\nu$ with
$|S_\nu|=1$. Put $S=\lim S_\nu$.
Let $I=[\Im S]=I(\hsig), \ K=[\Ker S]=K(\hsig) \subset \P^N$ 
be  the limit image and the limit kernel of $\{S_\nu\}$, respectively. 
If the algebraic set $I\cap \cG$ is of positive dimension,
then $B = \bigcup_{\hat \ell \in I\cap \cG}|\ell|$ 
is an algebraic manifold contained in $\Lambda(\Gamma)$
with dimension more than $n$.
This is absurd, since $\Omega(\Gamma)$ contains an $n$-plane
which  does not intersect $B$.
Hence $I \cap \cG$ is a finite set. 
Consequently, $I \cap \cG$ consists of a single point, 
since it is the set of limit points of $\cG \setminus K$, 
which is connected. 
\ \qed

\begin{prop}\label{LimitIsOnePt} 
The limit image $I(\hsig)$ consists of a single point in $\cG$ 
for any normal sequence $\ssig$ in $\Gamma$.
\end{prop}
{\bf Proof} \  
We use the notation in the proof of the lemma above. By the lemma, we have 
$I\cap\cG=\{\hat\ell\}$ for some point $\hat\ell\in\cG$. 
Suppose that $\dim I > 0$. The linear map $S$ defines the projection 
$S : \P^N-K \to I$.
Since $\cG$ is not contained in any proper linear subspace in $\P^N$,  
there is a point $w \in I\setminus \cG$. 
The fiber $S^{-1}(w)$ does not intersects $\cG$ outside $K$, 
since otherwise, for $x \in \cG\setminus K$, 
we have $w=S(x)=\lim_\nu S_\nu(x) \in \cG$. This is absurd.   
Thus $\cG \subset K\cup S^{-1}(\hat\ell)$. This contradicts again 
the fact that the manifold $\cG$ is not contained in any proper
linear subspace in $\P^{N}$. Hence
we have $I=I\cap\cG=\{\hat\ell\}$. Thus we have the proposition. 
\ \qed

\begin{thm}\label{MasD5} 
Let $\ssig$ be a sequence of distinct elements of $\Gamma$. 
Then there are limit $n$-planes $\ell_I$, $\ell_K$, 
and a subsequence $\{\tau_\nu\}$ of $\ssig$,  
such that $\{\tau_\nu\}$ is uniformly 
convergent to $\ell_I$ on $\P^{2n+1}\setminus\ell_K$ 
in the following sense that, 
for any compact subset $M \subset \P^{2n+1}\setminus \ell_K$, 
and, for any neighborhood $V$ of $\ell_I$,
there is an integer $m_0$ such that
$\tau_\nu(M) \subset V$ for any $m>m_0$.
\end{thm}
{\bf Proof} \  
Choose a normal subsequence $\{\tau_\nu\}$ of $\{\sigma_\nu\}$ such that
$\{\tau^{-1}_\nu\}$ also has a convergent sequence of representatives.
Let $\{T_\nu\} \subset \GL_{N+1}(\C)$ be the convergent sequence corresponding to 
$\{\hat\tau_\nu\}$. Put $T=\lim_\nu T_\nu$. 
Note that $\{T'_\nu\}$, $T'_\nu=|T^*_\nu|^{-1}T^*_\nu$, 
is a convergent sequence of represents of $\{\hat\tau^{-1}_\nu\}$ by 
$(\ref{Astar})$ and $(\ref{inverse})$. 
Hence $\{\tau^{-1}_\nu\}$ is also a normal sequence. 
Put $T'=\lim_{\nu}T'_\nu$. 
By Proposition \ref{LimitIsOnePt}, $[\Im T]$ is a single point in $\cG$, 
which is corresponding to a limit $n$-plane, denoted by $\ell_I$, in $\P^{2n+1}$. 
On the other hand, 
since $[\Im T']$ is the limit image of the normal sequence $\{\hat\tau^{-1}_\nu\}$, 
$[\Im T']$ consists of a single point corresponding to a limit $n$-plane 
in $\P^{2n+1}$ by Proposition \ref{LimitIsOnePt}, which we denote by  $\ell_K$.
Note that $[\Im T']^{\perp}$ is the set of points parameterizing 
$n$-planes intersecting $\ell_K$ by Proposition \ref{PQXY}.
Since $\Ker T=(\Im T')^\perp$ by Theorem \ref{IinQ}, and
$\{\hat\tau_\nu\}$ converges uniformly on compact sets in 
$\P^N\setminus [\Ker T]$ to $[\Im T]$ by Theorem \ref{Myr1}, 
we see that $\{\tau_\nu\}$ 
converges uniformly compact sets on $\P^{2n+1}\setminus\ell_K$ to $\ell_I$.
This proves the theorem.
\ \qed

\bigskip
In the course of the proof, we have shown the following.
\begin{prop}\label{MasD2} 
Let $\ell_0$ be a limit $n$-plane of $\Gamma$. 
Then there are a limit $n$-plane 
$\ell_\infty$, and a normal sequence $\ssig \subset \Gamma$ such that 
$\ssig$ is uniformly convergent to $\ell_0$ on any compact set 
in $\P^{2n+1}\setminus\ell_\infty$ and that 
$\{\sigma^{-1}_\nu\}$ is uniformly convergent to 
$\ell_\infty$ on any compact set in 
$\P^{2n+1}\setminus\ell_0$.
\end{prop}

Next we shall show the following.
\begin{thm}\label{MasD6} For a type $\ELL$ group $\Gamma$, 
$\Lambda(\Gamma)$ is a closed, 
nowhere dense  $\Gamma$-invariant subset in $\P^{2n+1}$.
\end{thm}
{\bf Proof} \ . 
To show that $\Lambda(\Gamma)$ is $\Gamma$-invariant, 
we take any point $x\in\Lambda(\Gamma)$. 
Since $x$ is on a limit $n$-plane, say $\ell_0$, 
there is a normal sequence $\ssig$ of $\Gamma$ with 
$I(\hsig)=\hat\ell_0$ by Proposition \ref{MasD2}.  
Then $\{\sigma\circ\sigma_\nu\}$ is a normal sequence with 
$I(\{\hat\sigma\circ\hat\sigma_\nu\}) = \hat\sigma(\hat\ell_0)$.
Since the limit $n$-plane $\sigma(\ell_0)$ passes through the point 
$\sigma(x)$, 
$\Lambda(\Gamma)$ is $\Gamma$-invariant. 

To show that $\Lambda(\Gamma)$ is closed, 
let $\{x_\nu\}$ be a sequence of points of $\Lambda(\Gamma)$ 
such that $\lim_\nu x_\nu = x$ for some point $x \in \P^{2n+1}$. 
Let $\ell_\nu$ be a limit $n$-plane through $x_\nu$. 
By Proposition \ref{MasD2}, for each $\nu$, 
we can find a limit $n$-plane $\ell_{\nu,\infty}$ and a normal sequence
$\{\sigma_{\nu,k}\}_k$ such that 
$I(\{\hat\sigma_{\nu,k}\}_k) = \hat\ell_\nu$, and
that the sequence $\{\sigma_{\nu,k}\}_k$ is uniformly
convergent to $\ell_\nu$ on compact sets in 
$\P^{2n+1}\setminus\ell_{\nu,\infty}$. 
Taking a subsequence of $\{\ell_\nu\}$,
we can assume that the $\ell_\nu$ are all distinct and that 
$\{\hat\ell_\nu\}$ and $\{\hat\ell_{\nu,\infty}\}_\nu$ are convergent 
in $\cG$. 

Since $\{\hat\ell_{\nu,\infty}\}_\nu$ is convergent, there is an $n$-plane $\ell_a$
which is disjoint from the closure of $\bigcup_\nu |\ell_{\nu,\infty}|$.
Take a small tubular neighborhood $W$ of $\ell_a$, which is biholomorphic to the 
domain $(\ref{U})$, such that the closure $[W]$ is still disjoint from 
the closure of $\bigcup_\nu |\ell_{\nu,\infty}|$.

Fix a metric on $\cG$ and consider distance of points on $\cG$. 
Let $\delta_\nu$ be the minimal distance from 
$\hat\ell_\nu$ to any other $\hat\ell_\mu$ in $\cG$. 
Obviously, $\lim_\nu\delta_\nu=0$. Set 
\[ N_{\delta_\nu}(\hat\ell_\nu) = \{ z \in \cG : \dist(z, \hat\ell_\nu) < \delta_\nu \}. \]
Choose $k(\nu)$ such that 
\[ \hat\sigma_{\nu,k(\nu)}([W])\subset N_{\delta_\nu}(\hat\ell_\nu), \]
and that the $\sigma_{\nu,k(\nu)}$, $\nu=1,2,3,\dots$, are all distinct.
Put $\hat\ell=\lim_\nu\hat\ell_\nu$. Take any $\delta>0$. 
Then there is $\nu_0$ such that 
$N_{\delta_\nu}(\hat\ell_\nu) \subset N_\delta(\hat\ell)$ holds 
for any $\nu>\nu_0$.
Thus, for $\nu>\nu_0$, we have 
$\hat\sigma_{\nu,k(\nu)}([W])\subset N_{\delta}(\hat\ell)$. 
This implies that 
$\{\hat\sigma_{\nu,k(\nu)}\}$ converges to $\hat\ell$ uniformly on $W$. 
Thus, $\ell$ is a limit $n$-plane passing through $x$.
Hence $\Lambda(\Gamma)$ is closed. 

Lastly, we shall show that $\Lambda(\Gamma)$ is nowhere dense. 
Let $x$ be any point in $\Lambda(\Gamma)$. 
By Proposition \ref{MasD2}, there are $n$-planes 
$\ell_0$, $\ell_\infty$ in $\P^{2n+1}$ and a normal sequence $\{\sigma_\nu\}$ 
such that $x \in \ell_0$ and that $\lim_\nu \hat\sigma_\nu(\hat K) = \hat\ell_0$ for any 
compact set $K \subset \P^{2n+1}\setminus \ell_\infty$. 
By the property $\ELL$, we can set $K$ as a  single $n$-plane $\ell$ contained  in $\Omega(\Gamma)$.
Then, for every neighborhood $W$ of $x$, 
there is an integer $\nu_0$ such that $W \cap \sigma_\nu(\ell)\neq \emptyset$ for $\nu\geq \nu_0$. 
Thus $W$ contains a point in $\Omega(\Gamma)$. Hence 
$\Lambda(\Gamma)$ is nowhere dense. 
\ \qed

\begin{thm}\label{CanQ} For a type $\ELL$ group $\Gamma$, 
the action of $\Gamma$ on $\Omega(\Gamma)$ is properly discontinuous.
\end{thm}
{\bf Proof} \  Take any compact set $M$ in $\Omega(\Gamma)$. 
Suppose that there is 
an infinite sequence $\{\sigma_\nu\}_\nu$ of distinct elements of $\Gamma$ 
such that $M \cap \sigma_\nu(M) \neq \emptyset$ for any $\nu$. 
By Proposition \ref{MasD2}, 
replacing $\ssig$ with its normal subsequence, we can assume that 
there are limit $n$-planes $\ell_K$ and $\ell_I$ such that $\ssig$ converges uniformly on $\P^{2n+1}\setminus\ell_K$
to $\ell_I$.  Since $\Omega(\Gamma)$ has no intersection with limit 
$n$-planes, we see that 
$M\cap(\ell_I \cup \ell_K)=\emptyset$. Therefore $\{\sigma_\nu(M)\}$ converges to a subset on $\ell_I$.
This contradicts the assumption that $M \cap \sigma_\nu(M) \neq \emptyset$ for any $\nu$. 
\ \qed

\bigskip
By Theorem \ref{CanQ}, we can define canonically the quotient space 
$\Omega(\Gamma)/\Gamma$, which we denote by $X(\Gamma)$,
\[ X(\Gamma)= \Omega(\Gamma)/\Gamma. \]

\begin{rem}
There are examples of $\Gamma$ for which $X(\Gamma)$ is not connected. 
Such an example can be constructed easily in case $n=1$ by considering 
a flat twistor space over a conformally flat 4-manifold (\cite{ext}), where 
every connected component of $\Omega(\Gamma)$ is large. 
We do not know, however, whether this is the case for all type $\ELL$ groups 
or not.
\end{rem}

\section{Discontinuous group actions on large domains}
\label{setting1}
In this section, we shall show that a large domain which covers a compact manifold is a connected component of
$\Omega(\Gamma)$ of some $\Gamma$ of type $\ELL$. 

\begin{prop}\label{largeL} 
Let $\Gamma$ be a group of holomorphic automorphisms 
of a large domain $\Omega$ in $\P^{2n+1}$. 
Suppose that $\Gamma$ is torsion free and that 
the action of $\Gamma$ on $\Omega$   
is properly discontinuous. 
Then $\Gamma$ is of type $\ELL$.
\end{prop}
{\bf Proof} \  
First we shall prove that $\Gamma$ is a subgroup of $\PGL_{2n+2}(\C)$.
By a {\em line} we shall mean 
a $1$-dimensional projective linear subspace of a projective space.
Since every line in $\P^{2n+1}$ has a tubular neighborhood 
with a smooth convex-concave boundary, 
the following lemma follows immediately 
from a theorem of Ivashkovich \cite{Iv}.
\begin{lem} \label{ext2PGL} 
Let $L_\nu$, $\nu=1,2$, be lines in $\P^m$ $(m\geq 2)$, and 
$U_\nu$ a tubular neighborhood of $L_\nu$. 
Suppose that $\gamma : U_1 \to U_2$
is a biholomorphic mapping. Then $\gamma$ extends to an element of 
$\PGL_{m+1}(\C)$.
\end{lem}

\begin{lem} \label{rank} 
Let $\sigma\in\Gamma\setminus\{1\}$ be any element, 
and $\tsigma \in \GL_{2n+2}(\C)$ a representative of $\sigma$. Then, the inequality
\begin{equation}\label{nn} 
 \rank(\tsigma-\al I) \geq n+1
\end{equation}
holds for any $\al \in \C$.
\end{lem}
{\bf Proof} \   
Consider the subspace  
\begin{equation}\label{V} 
V=\{ z \in \C^{2n+2} : (\tsigma - \al I)z=0 \}. 
\end{equation}
Each point of the projectivised linear subspace $[V]\subset \P^{2n+1}$ 
is fixed by $\sigma$. Suppose that $\rank(\tsigma-\al I)\leq n$.
Then $\dim V \geq n+2$. Therefore any $n$-plane in $\Omega$ intersects $[V]$
and the every point on the intersection is fixed by $\sigma$.
This is absurd, since $\Gamma$ is torsion free and properly discontinuous 
on $\Omega$. Thus we have the lemma. 
\ \qed

\begin{lem} \label{IIL} \footnote{Compare \cite[Lemma 1.6]{Lar}, which is for
$1$-planes in $\P^m$} If $(\ref{nn})$ holds for
any $\al \in \C$, then there is an $n$-plane $\ell$ such 
that $\sigma(\ell) \cap \ell = \emptyset$.
\end{lem}
{\bf Proof} \   We have to choose a 
subspace $L\subset \C^{2n+2}$ of dimension $n+1$ such that 
$\tsigma(L) \cap L =\{0\}$.
Put  
\[ \rho = \min_{\al \in \C}\rank(\tsigma - \al I). \]
We can assume that $\rho$ is attained at $\al=1$ without loss of generality.
We put $N=\tsigma - I$ and then $\rho=\rank N$.  
Define $\vf : C^{2n+2}\to \C^{2n+2}$ by $\vf(z)=Nz$. 
Since $\rho\geq n+1$ by the assumption, there is an $(n+1)$-dimensional subspace 
$L_1\subset \Im\vf$. Put $\tilde{L}_1=\vf^{-1}(L_1)$.
Then, since $\dim \Ker\vf = 2n+2 - \rho$, we have $\dim\tilde{L}_1=3n+3-\rho$. 
Since $\dim L_1= n+1$ and $\dim\Ker\vf=2n+2-\rho$,
we can choose a subspace $L \subset \tilde{L}_1$ 
such that $\dim L= n+1$, $L\cap L_1 = \{0\}$ and $L\cap\Ker\vf=\{0\}$.
We claim that $L$ is the desired linear subspace in $\C^{2n+2}$. 
To verify the claim, we choose $X \in M((2n+2)\times (n+1),\C)$ with $\rank X=n+1$ 
such that
\[ L = \{ z \in \C^{2n+2} : z=Xu, \ \ u \in \C^{n+1}\}. \]
Then $L\cap \tsigma(L)=\{0\}$ holds if and only if
\[
\det (\tsigma X, X) \neq 0. 
\]
This is equivalent to 
\[
\det (NX,X) \neq 0. 
\]
That $L\cap\Ker\vf=\{0\}$ implies that $NX$ is of maximal rank, and  
that $L\cap L_1 = \{0\}$ implies that the vectors in $NX$ and $X$ span $\C^{2n+2}$.
Thus the claim is verified. 
\ \qed

\bigskip
Now we go back to the proof of Proposition \ref{largeL}. 
By the assumption that $\Omega$ is large, 
there is a relatively compact subdomain $W \subset \Omega$ 
which is biholomorphic to $U$. The $n$-planes in $W$ are parametrized by $\hat W \subset \cG\subset \P^N$. 
Since the action of $\Gamma$ on $\Omega$ is properly discontinuous, the set
\[ S = \{ \sigma \in \Gamma\setminus \{1\} : \hat\sigma(\hat W) \cap \hat W \neq \emptyset \} \]
is finite. 
Let $\ell$ be an $n$-plane in $W$. 
For $\sigma \in S$,  we have $Q(\hat\ell, \hat\sigma(\hat\ell))=0$,
when $\ell$ intersects $\sigma(\ell)$.  
By Lemmas \ref{rank} and \ref{IIL}, we see that the set
\[ Y_\sigma = \{ \zeta \in \cG : Q(\zeta, \hat\sigma(\zeta))=0 \}  \]
is a proper analytic subset of $\cG$. Hence the set
\[ V= \hat W \setminus \bigcup_{\sigma \in S} Y_\sigma \]
is not empty. Take a point $\hat \ell' \in V$. 
Then, we can choose a neighborhood $W'$ of $\ell'$ which is biholomorphic to $U$ and
satisfies $\sigma(W')\cap W'=\emptyset$ for all $\sigma$ in $S$, and hence in $\Gamma$.  
\ \qed

\begin{thm}\label{component}
Let $\Omega \subset \P^{2n+1}$ be a large domain which
is an unramified cover of a compact complex manifold. 
Then there is a type $\ELL$ group $\Gamma$ such that 
$\Omega(\Gamma)$ contains $\Omega$ as a connected component.
\end{thm}
{\bf Proof} \   By Lemma \ref{ext2PGL}, there is a group 
$\Gamma\subset \PGL_{2n+2}(\C)$ of 
holomorphic automorphisms of $\Omega$ such that $\Omega/\Gamma$ 
is compact. Since $\Gamma$ is finitely generated, we 
can assume that $\Gamma$ is torsion free by Selberg's lemma. 
Hence by Proposition \ref{largeL}, $\Gamma$ is of type $\ELL$.

We claim that $\Omega \subset \Omega(\Gamma)$.
To verify this, suppose contrarily that there is a point 
$x \in \Omega \cap \Lambda(\Gamma)$. 
Then there are limit $n$-planes $\ell_I, \ell_K$ such that $x\in \ell_I$,
and a sequence $\{\sigma_m\}$ of distinct elements of $\Gamma$ such that 
$\{\sigma_m\}$ converges uniformly on $\P^{2n+1}\setminus\ell_K$ to $\ell_I$. 
Let $\ell$ be an $n$-plane contained in $\Omega$. 
Displacing $\ell$ a little if necessary, we can assume that $\ell\cap\ell_K=\emptyset$.
Let $K_x$ be a compact neighborhood of $x$ contained in $\Omega$. 
Put $K=K_x\cup\ell$, which 
is a compact set contained in $\Omega$. Since $\{\sigma_m(\ell)\}$ 
converges to $\ell_I$,
we see that $\sigma_m(K) \cap K \neq \emptyset$ for infinitely many $m$. 
This contradicts the assumption that $\Gamma$ is properly discontinuous on $\Omega$. 
Thus the claim is verified.

Now $\Omega$ is contained in a connected component, say $\Omega_0$, 
of $\Omega(\Gamma)$.  
Since $\Omega$ is $\Gamma$-invariant, so is $\Omega_0$. 
Therefore, by Theorem \ref{CanQ},  
$\Omega_0/\Gamma$ is a connected complex spaces which contains $\OG$.
Since $\OG$ is compact, we infer that $\OG=\Omega_0/\Gamma$. 
Hence $\Omega=\Omega_0$. 
\ \qed

\begin{prop}\label{Kaehler}
Let $X$ be a compact K\"ahler manifold which contains 
a domain $W$ biholomorphic to 
\[ 
U=\{ [z_0:\dots:z_{2n+1}] \in \P^{2n+1}
 : |z_0|^2 +\cdots + |z_n|^2 < |z_{n+1}|^2 +\cdots + |z_{2n+1}|^2 \}
\]
Then $X$ is unirational. In particular, $X$ is simply connected.
\end{prop}
{\bf Proof} \   The proof of \cite[Corollary 3.1]{line} works also in this case.  
Take any $n$-plane $\ell \subset W$. 
Let $B$ be the irreducible component of the Barlet space which contains the 
point $\hat\ell$ corresponding to $\ell$. Since $X$ is K\"ahler, $B$ is compact.
Consider the graph
\[ Z = \{ (x, b) \in X\times B \ : \ x \in b \} \]
Let $p_X : Z \to X$ and $p_B : Z \to B$ the natural projections. Fix a point $o\in\ell$. 
Since $B$ is compact, 
we can apply a theorem of Campana \cite[Corollaire 1]{Camp}, 
which says 
that $p_X^{-1}(o)$ is a compact {\it algebraic} variety.
Hence $B_o:=p_B(p_X^{-1}(o))$ is also compact and algebraic.
Put $M:=p_B^{-1}(B_o)$ and $f=P_B|_M$.
Then $f : M \to B_o$ is a $\P^n$-fiber space over a compact algebraic variety $B_o$.
By the choice of $o$, $B_o$ is non-singular at $o$, 
and there is a small open neighborhood $N\subset B_o$ centered at 
$\hat\ell$, and a biholomorphic map $\tau : f^{-1}(N) \to N\times \P^n$ such that 
$f = p\circ \tau$, where $p : N\times \P^n\to N$ is the projection.
Let $\mu : M^*\to M$ be a desingularization of $M$ which is a succession of blowing-ups.  
Here we can assume that $\mu$ is biholomorphic on $f^{-1}(N)$. 
Thus we have a fiber space $g:=f\circ \mu : M^* \to B_o$ whose general fiber is $\P^n$.  

\begin{lem}\label{Malg}
$M^*$ is an algebraic variety.
\end{lem}
{\bf Proof} \  
As in \cite[\S 12]{U}, we consider the direct image sheaf $g_*\cO(-K_{M^*})$, 
and the associated projective fiber space $\P(g_*\cO(-K_{M^*}))$ over $B_o$.
Note that $\P(g_*\cO(-K_{M^*}))$ is an algebraic space. 
Since $g_*\cO(-K_{M^*})$ is a locally free sheaf of $\rank={_{2n+1}C_{n+1}}$ 
on a non-empty Zariski open subset of $B_o$, we have a commutative diagram
\[
\begin{array}{ccc}
M^* \ \ \ & \stackrel{h}{\longrightarrow} & \P(g_*\cO(-K_{M^*})) \\
g \searrow & & \swarrow \pi \\
& B_o, & 
\end{array}
\]
where $h$ is a meromorphic map whose restriction
\[ h|_{g^{-1}(b)} : g^{-1}(b) \to \P(g_*\cO(-K_{M^*}))_b \]
to a general fiber $g^{-1}(b)$ is the map defined by the linear system $|\cO_{\P^n}(n+1)|$.
Hence we infer that $\dim h(M^*)=\dim M^*$. Since
$\P(g_*\cO(-K_{M^*}))$ is algebraic, so is $h(M^*)$. Hence $M^*$ is algebraic. 
\ \qed

Since $X=p_X(M)=p_X(\mu(M^*))$, we see that $X$ is algebraic by Lemma \ref{Malg}. 
Let $\nu: Y \to X$ be a succession of blowing-ups such that $Y$ is projective algebraic. 
Let $j : U \to W$ be a biholomorphic map.
Since any meromorphic function on $U$ extends to a meromorphic function on 
$\P^{2n+1}$, $\nu^{-1}\circ j : U \to Y$ extends to a meromorphic map 
$\P^{2n+1} \cdots\!\!> Y$. This implies that $Y$ is unirational. Hence $X$ is unirational. 
\ \qed

\begin{thm}\label{nonKaehler}
A compact complex manifold which is covered by  
a large domain in $\P^{2n+1}$ is non-K\"ahler, except for $\P^{2n+1}$ itself.
\end{thm}
{\bf Proof} \  This follows from Proposition \ref{largeL}, Proposition \ref{Kaehler} 
and Theorem \ref{component}. 
\ \qed

\bigskip
Note that, for a large domain in Theorem \ref{nonKaehler}, 
we assume nothing on its fundamental group
nor on its complement in $\P^{2n+1}$. 
Thus our result gives a slight generalization of  \cite[Proposition 1.9]{Lar} 
for odd dimensional projective spaces.

\section{Klein combinations}
Let $\Omega_\nu \subset \P^{2n+1}$, $\nu=1,2$, be large domains, and
$\Gamma_\nu \subset \Aut(\Omega_\nu)$ free and properly discontinuous groups.
Put
\[ 
U(\ep) = \{ |z_0|^2 +\cdots+|z_n|^2 < \ep(|z_{n+1}|^2 +\cdots+ |z_{2n+1}|^2) \} \subset \P^{2n+1}, 
\ \ \ \ep>1,
\]
\[ 
N(\ep) = \left[U(\ep)\right] \setminus U(\ep^{-1}). 
\]
Then 
\[ 
\sigma : N(\ep) \to N(\ep), \ \ 
\sigma([z_0:\dots:z_n:z_{n+1}:\dots:z_{2n+1}])
=[z_{n+1}:\dots:z_{2n+1}:z_0:\dots:z_n]
\]
is a biholomorphic map. 
Let $j_\nu : U(\ep) \to X_\nu=\Omega_\nu/\Gamma_\nu$  
be holomorphic open embeddings.
Then we can consider the gluing 
\[
X_1\#X_2 = \left(X_1\setminus j_1(U(\ep^{-1}))\right) \bigcup \left( X_2\setminus j_2(U(\ep^{-1})) \right)
\]
by $ j_2\circ \sigma \circ {j_1}^{-1} : j_1(N(\ep)) \to j_2(N(\ep))$ 
to obtain a new complex manifold. Then 
\[  X_1 \# X_2 = \Omega/\Gamma \]
for some large domain $\Omega \subset \P^{2n+1}$ and $\Gamma$ (\cite{line}).
Here we have $\Gamma \simeq \Gamma_1\ast \Gamma_2$. 
$X_1\# X_2$ is called the {\em Klein combination} of $X_1$ and $X_2$. 
If $\Gamma_\nu$'s are cocompact then
so is $\Gamma$ on $\Omega$. 

\bigskip
The {\em handle attachments} can also be defined. In those cases, we have 
$\Gamma \simeq \Gamma_1\ast \Z$. 
Thus we can get many examples of $\Gamma$ and $X(\Gamma)$.

\section{An analogue of the Ford region}
Fix a system of homogeneous coordinates 
$[z^0:z^1:\dots:z^n:z^{n+1}:\dots:z^{2n+1}]$ on $\P^{2n+1}$.
Put $z'=(z^0,\dots,z^n), \ z''=(z^{n+1},\dots,z^{2n+1})$, and
write $[z':z'']$ instead of 
$[z^0:z^1:\dots:z^n:z^{n+1}:\dots:z^{2n+1}]$ for brevity.
Let $\ell''$ be the $n$-plane defined by $z''=0$.  
Put 
\[ E=\P^{2n+1}\setminus \ell'', \]
and define the projection by
\[ \pi : E \to \P^n, \ \ \pi([z':z''])=z''. \]
Then $E$ is isomorphic to 
$\cO_{\P^n}(1)^{\oplus(n+1)}$ as a vector bundle over $\P^n$. 

\subsection{Volume form on $E$}
Take the open covering of $E=\bigcup_{\al=1}^{n+1}U_\al$, where   
\[ U_\al=\{ [z':z'']\in \P^{2n+1} : z^{n+\al}\neq 0 \}, \ \ 1\leq\al\leq n+1. \]
On each $U_\al$, we define a system of coordinates by  
\[ \left\{
\begin{array}{lcll}
\zeta_\al^j & = & \Frac{z^j}{z^{n+\al}}, & 0\leq j\leq n, \\
\\
x_\al^k & = & \Frac{z^{n+k}}{z^{n+\al}}, & 1\leq k < \al,\\
\\
x_\al^{k-1} & = & \Frac{z^{n+k}}{z^{n+\al}}, & \al < k \leq n+1.
\end{array}\right.
\]
Then $\pi|U_\al$ is given by
\[ 
\pi(\zeta_\al^0,\dots, \zeta_\al^n , x_\al^1,\dots, x_\al^n)
=(x_\al^1,\dots, x_\al^n).
\]
On $U_\al$, we define
\begin{eqnarray}
 d\zeta_\al & = & d\zeta_\al^0\wg\dots\wg d\zeta_\al^n\\ 
 dx_\al & = & dx_\al^1\wg\dots\wg dx_\al^n 
\end{eqnarray}
and put
\begin{equation}
dV_\al =\sqrt{-1}(1+\|x_\al\|^2)^{-2(n+1)}
d\zeta_\al\wg \overline{d\zeta_\al}\wg dx_\al\wg\overline{dx_\al},
\end{equation}
where
\[ \|x_\al\|^2 = \sum_{k=1}^n|x_\al^k|^2. \]

It is easy to check that the $(2n+1,2n+1)$-forms $dV_\al$ 
patch together to give a global volume form 
\begin{equation}
dV = dV_\al \mbox{  on  } U_\al
\end{equation}
on $E=\P^{2n+1}\setminus \ell''$.

\begin{lem}\label{proj} 
Consider the projective transformation of $\P^m$ defined by 
\[ 
y^\lm=\Frac{c^\lm_\mu x^\mu+c_0^\lm}{c_\mu x^\mu+c_0},
\ \  1\leq \lm, \mu \leq m.
\]
Then we have
\[ dy^1\wg\dots\wg dy^m
=\frac{\det C}{(c_\mu x^\mu+c_0)^{m+1}}
dx^1\wg\dots\wg dx^m,
\]
where 
\[ 
C=\left(\begin{array}{cccc}
c_0        & c_1 & \dots & c_m \\
c_0^1 & c^1_1 & \dots & c^1_m \\ 
\vdots &          & \vdots \\
c_0^m & c^m_1 & \dots & c^m_m
\end{array}\right)
\]
\end{lem}
{\bf Proof} \  Put 
\[ P=c_\mu x^\mu+c_0, \ \ Q=c_\mu x^\mu, \ \ 
p^\lm=c^\lm_\mu x^\mu+c^\lm_0, \ \ 
q^\lm=c^\lm_\mu x^\mu, 
\]
where $\mu$ is summed for $\mu=1,\dots,m$. Then we have
\begin{eqnarray*}
dy^1\wg\dots\wg dy^m 
& = & 
\Wg^m_{\lm=1}\left(P^{-1}dq^\lm - p^\lm P^{-2}dQ\right)
\\  & = & 
P^{-2m}\Wg^m_{\lm=1}\left(Pdq^\lm - p^\lm dQ\right)
\\  & = &
P^{-(m+1)}\left(
P\Wg^m_{\lm=1}dq^\lm 
+ \sum_{k=1}^m
(-1)^kp^k dQ\wg dq^1\wg\dots\wg dq^{k-1}\wg
dq^{k+1}\wg\dots\wg dq^m\right). 
\end{eqnarray*} 
Define $A$ and $A_k$ by 
\begin{eqnarray*} 
A dx^1\wg\dots\wg dx^m & = &dq^1\wg\dots\wg dq^\lm,
\\ 
A_k dx^1\wg\dots\wg dx^m & = & dQ\wg dq^1\wg\dots\wg dq^{k-1}\wg dq^{k+1}\wg\dots\wg dq^m.  
\end{eqnarray*}
Then we have
\begin{equation}\label{dyAdx}
dy^1\wg\dots\wg dy^m=
P^{-(m+1)}\left((c_\mu x^\mu+c_0)A+\sum_{k=1}^m(-1)^k(c^k_\mu x^\mu+c^k_0)A_k\right)dx^1\wg\dots\wg dx^m\end{equation}
Note that 
\[
A=\det \left(\begin{array}{ccc}
c^1_1 & \dots & c^1_m \\ 
\vdots &          & \vdots \\
c^m_1 & \dots & c^m_m
\end{array}\right)
\mbox{  and }
A_k=\det \left(\begin{array}{ccc}
c_1 & \dots & c_m \\
c^1_1 & \dots & c^1_m \\ 
\vdots &          & \vdots \\
c^{k-1}_1 &          & c^{k-1}_m \\
c^{k+1}_1 &          & c^{k+1}_m \\
\vdots &          & \vdots \\
c^m_1 & \dots & c^m_m 
\end{array}\right)
\]
Thus, we have 
\[ 
c_\mu A + \sum_{k=1}^m(-1)^kc_\mu^kA_k =
\det
\left(\begin{array}{cccc}
c_\mu        & c_1 & \dots & c_m \\
c_\mu^1 & c^1_1 & \dots & c^1_m \\ 
\vdots &          & \vdots \\
c_\mu^m & c^m_1 & \dots & c^m_m
\end{array}\right)=0
\]
for $\mu=1,\dots, m$, and 
\[ 
c_0 A + \sum_{k=1}^m(-1)^kc_0^kA_k =
\det
\left(\begin{array}{cccc}
c_0        & c_1 & \dots & c_m \\
c_0^1 & c^1_1 & \dots & c^1_m \\ 
\vdots &          & \vdots \\
c_0^m & c^m_1 & \dots & c^m_m
\end{array}\right)=\det C.
\]
Hence it follows from (\ref{dyAdx}) that
\[
dy^1\wg\dots dy^m
=P^{-(m+1)}\det C dx^1\wg\dots\wg dx^m. 
\]
\ \qed

\bigskip
\begin{lem}\label{mu} 
\footnote{This is the corrected version of \cite[Lemma 3.2]{Leb}. 
There was a mistake in the calculation there. 
The results \cite[Proposition 3.1, Lemma 3.3]{Leb} 
hold true. Calculations in the proofs there should be corrected accordingly, 
but need no essential changes.  Sublemmas $3.1$,\ $3.2$ in \cite{Leb} 
and their proofs are correct.}
For 
\[ g=\MAT{A}{B}{C}{D} \in \SL_{2n+2}(\C), 
\ \ \ A, B, C, D \in M_{n+1}(\C)
\]
with $\det C\neq 0$, the pull-back of $dV$ is given by
\[ g^*dV=\mu^{4(n+1)}_g dV \]
where
\[
\mu_g(z)=\Frac{\|z''\|}{\|Cz'+Dz''\|}.
\]
\end{lem}
{\bf Proof} \  
We write a square matrix $M$ of size $(n+1)$ as 
\[ 
M=\left(\begin{array}{ccc}
m^0_0  & \dots & m^0_n \\
\vdots &          & \vdots \\
m^n_0 & \dots & m^n_n 
\end{array}\right).
\]
Set $\al=n+1$ and consider the projective transformation $g$ on 
$U_\al=U_{n+1}$. We omit the subscript $n+1$ for simplicity, 
and write the local coordinates by $(\zeta^0,\dots,\zeta^n, x^1,\dots, x^n)$ 
instead of 
$(\zeta^0_{n+1},\dots,\zeta^n_{n+1}, x^1_{n+1},\dots, x^n_{n+1})$.
Then  $g$ sends $(\zeta^j, x^k)$ to $({\zeta^j}',{x^k}')$, where 
\begin{eqnarray*}
 {\zeta^j}' & = & \Frac{\sum_{\lm=0}^n a^j_\lm \zeta^\lm + \sum_{\mu=0}^{n-1} b^j_\mu x^{\mu+1}+b^j_n }
{\sum_{\lm=0}^n c^n_\lm \zeta^\lm + \sum_{\mu=0}^{n-1} d^n_\mu x^{\mu+1}+d^n_n}, \ \ j=0,\dots,n, \\
& & \\
 {x^k}' & = & \Frac{\sum_{\lm=0}^n c^k_\lm \zeta^\lm + \sum_{\mu=0}^{n-1} d^k_\mu x^{\mu+1}+d^k_n }
{\sum_{\lm=0}^n c^n_\lm \zeta^\lm + \sum_{\mu=0}^{n-1} d^n_\mu x^{\mu+1}+d^n_n}, \ \ k=1,\dots,n. 
\end{eqnarray*}

Then by Lemma \ref{proj}, we have
\begin{equation}
d{\zeta}'\wg\overline{d\zeta'}\wg dx'\wg\overline{dx'}=
\left|\sum_{\lm=0}^n c^n_\lm \zeta^\lm + \sum_{\mu=0}^{n-1} d^n_\mu x^{\mu+1}+d^n_n\right|^{-4(n+1)} 
d{\zeta}\wg\overline{d\zeta}\wg dx\wg\overline{dx}
\end{equation}
Hence we have
\begin{eqnarray*}
g^*dV & = &  \sqrt{-1}
\left( 1+ \sum_{k=1}^n\left|\Frac{\sum_{\lm=0}^n c^k_\lm \zeta^\lm + \sum_{\mu=0}^{n-1} d^k_\mu x^{\mu+1}+d^k_n }
{\sum_{\lm=0}^n c^n_\lm \zeta^\lm + \sum_{\mu=0}^{n-1} d^n_\mu x^{\mu+1}+d^n_n} \right|^2\right)^{-2(n+1)} \\
& & \times
\left|\sum_{\lm=0}^n c^n_\lm \zeta^\lm + \sum_{\mu=0}^{n-1} d^n_\mu x^{\mu+1}+d^n_n\right|^{-4(n+1)}\!\!
d{\zeta}\wg\overline{d\zeta}\wg dx\wg\overline{dx}\\
& = &  \sqrt{-1}
\left(\sum_{k=0}^n\left|\sum_{\lm=0}^n c^n_\lm \zeta^\lm + \sum_{\mu=0}^{n-1} d^n_\mu x^{\mu+1}+d^n_n\right|^2\right)^{-2(n+1)}\!\!
d{\zeta}\wg\overline{d\zeta}\wg dx\wg\overline{dx}\\
& =&  \sqrt{-1}
\|C\zeta + D\tx\|^{-4(n+1)}
d{\zeta}\wg\overline{d\zeta}\wg dx\wg\overline{dx},
\end{eqnarray*}
where $\tx=(x^1,\dots,x^n,1)$. Thus we have
\[
g^*dV =  
\left(\Frac{\|\tx\|}{\|C\zeta + D\tx\|}\right)^{4(n+1)}dV
=
 \left(\Frac{\|z''\|}{\|Cz' + Dz''\|}\right)^{4(n+1)}dV.
\]
This proves the lemma. 
\ \qed

\subsection{$F$-region}
Recall that the norm of $u=(u_1,\dots,u_m) \in \C^m$ is defined by
\begin{equation}\label{norm1}
\|u\|=\left(|u_1|^2+\dots+|u_m|^2\right)^{1/2}.
\end{equation}
The norm of a matrix $A=(a_{ij}) \in \Mt{m}$ is defined by the 
operator norm
\begin{equation}\label{norm2} 
\|A\| 
=\sup_{z\neq 0, z\in\C^m}\Frac{\|Az\|}{\|z\|}.
\end{equation}

Let $\Gamma \subset \PGL_{2n+2}(\C)$ be a type $\ELL$ group, 
and set $\Omega=\Omega(\Gamma)$, $\Lambda=\Lambda(\Gamma)$.
Put $\Gamma^*=\Gamma\setminus\{1\}$.
Recall the proof of Proposition \ref{largeL}, where it is shown that $Y_\sigma$ is 
proper analytic subset of $\cG$. That proof shows that, 
moving the $n$-plane $l''=\{z''=0\}$ slightly if necessary,  
we can choose a positive number $R$ such that the set 
\[  
V_R = \{ [z':z'']\in \P^{2n+1} : \|z'\| > R\|z''\| \}
\]
is contained in $\Omega$, and that 
\begin{equation}\label{gV_V}
  g(V_R) \cap V_R = \emptyset
\end{equation}
holds for any $g \in \Gamma^*$.
Every $g \in \Gamma$ has a representative $\tg \in \SL_{2n+2}(\C)$, 
which we write as 
\[ 
\tg = \MAT{A_g}{B_g}{C_g}{D_g}, \ \ \ A_g,  B_g, C_g, D_g \in 
\Mt{n+1}. 
\]

\begin{lem}\label{BdABCD}
There is a constant $R_0>0$ such that, for any $g \in \Gamma^*$, 
we have 

{\rm (i)} \ $\det C_g \neq 0$,\ \ \ 
{\rm (ii)} \ $\| A_gC_g^{-1} \| \leq R_0$,\ \ \ 
{\rm (iii)} \ $\| C_g^{-1}D_g \| \leq R_0$.
\end{lem}
{\bf Proof} \   We fix an $R_0=R$ which satisfies (\ref{gV_V}). 
${\rm (i)}$ Suppose that $\det C_g = 0$.  Then there
is a point $z$ on $l''$ such that $g(z) \in l''$. 
Thus $g(l'') \cap l'' \neq \emptyset$.
Since $g \neq 1$ by assumption, this contradicts (\ref{gV_V}).  
${\rm (ii)}$ The $n$-plane $g(l'')$ is 
given by $z'=A_gC_g^{-1}z''$.
Since $g(l'') \cap V_{R_0} = \emptyset$ by (\ref{gV_V}), 
we have $\|A_gC_g^{-1}\| \leq R_0$. 
${\rm (iii)}$ The equation of the $n$plane $g^{-1}(l'')$ is 
given by $z'= -C_g^{-1}D_gz''$, we 
have $\|C_g^{-1}D_g\| \leq R_0$ by the argument above. 
\ \qed

\bigskip
Put 
\[ \Delta_g = \{ z=[z':z''] \in \P^{2n+1} : \|z''\|< \|C_gz'+D_gz''\| \}, \]
\[ \bar{\Delta}_g = \{ z=[z':z''] \in \P^{2n+1} : \|z''\|\leq \|C_gz'+D_gz''\| \}, 
\]
\[ \Sigma_g = \{ z=[z':z''] \in \P^{2n+1} : \|z''\|= \|C_gz'+D_gz''\| \}, \]
\[ \Delta^c_g = \{ z=[z':z''] \in \P^{2n+1} : \|z''\|\geq \|C_gz'+D_gz''\| \}, \]
and 
\[ \bar{\Delta} = \bigcap_{g \in \Gamma^*} \bar{\Delta}_g. \]

\begin{definition}\label{Ford}
Consider the set of interior points 
\[ 
F = \Int \bar{\Delta}
\]
of $\bar\Delta$, which we call the $F$-region of the type 
$\ELL$ groups. 
\end{definition}
This is an analogue of the Ford region in the Kleinian group theory. 
In deed, for some type $\ELL$ groups which satisfy an additional condition 
(see $(\clubsuit), (\spadesuit)$ below), $F$ will give a fundamental set of the 
action of $\Gamma$ on $\Omega$. 
Now we put
\begin{eqnarray}
\label{truFord}  \bF & = & \mbox{the closure of $F$ in $\P^{2n+1}$,} \\
\nonumber \pd \bF & = & \bar F \setminus F. 
\end{eqnarray}

\bigskip
We consider the following set of positive real numbers:   
\begin{equation}\label{Cminus}
\cR=\{\|C_g^{-1}\| : g \in \Gamma^* \},
\end{equation}
and consider the conditions on $\cR$ :
\begin{enumerate}
\item[$(\clubsuit)$ ] $\cR$ is bounded in $\R$,
\item[$(\spadesuit)$ ] $\cR$ has 
no accumulation points other than $0$ in $\R$.
\end{enumerate}

\begin{rem}\label{invariant}
The number $\|C_g^{-1}\|$ is something like the radius of 
the isometric circle of $g$ in Kleinian group theory.
The conditions $(\clubsuit)$ and $(\spadesuit)$ may 
depend on the choice of homogeneous coordinates on $\P^{2n+1}$. 
But they are preserved under the coordinate change 
$w=\tau(z)$ of the form $\tau=\MAT{P}{Q}{0}{S}\in\PGL_{2n+2}(\C)$. 
\end{rem}

\begin{prop}\label{NbdOfL2}
The condition $(\clubsuit)$ is satisfied, if and only if 
$\bar F$ contains $V_R$ for some $R>0$. 
\end{prop}
{\bf Proof} \  To prove that $(\clubsuit)$ is sufficient, 
let $\rho>0$ be an upper bound of $\cR$.
Set $R=R_0+\rho$, 
where $R_0$ is the constant in Lemma \ref{BdABCD}(iii). 
Then, for any point $z=[z':z'']\in V_R$ and any $g\in\Gamma^*$, 
we have 
\[
\|C_gz'+D_gz''\|\geq 
\Frac{\|z'+C_g^{-1}D_gz''\|}{\|C_g^{-1}\|}
\geq 
\Frac{\|z'\|-R_0\|z''\|}{\|C_g^{-1}\|}
\geq 
\Frac{\rho\|z''\|}{\|C_g^{-1}\|}
\geq \|z''\|.
\]

To prove that $(\clubsuit)$ is necessary, take any $n$-plane 
\[ \ell_Y \ : \ z''=Yz', \ \ \ Y\in \Mt{n+1}, \  \|Y\|< R^{-1} \]
in $V_R$. Since $\ell_Y \subset \bDelta_g$ for any $g\in\Gamma^*$, we
have
\begin{equation}\label{YCD}
 \| Yz'\|\leq \|C_gz'+D_gYz'\| 
\end{equation}
for any $z'\in\C^2$.  Put 
\[ 
G = (C_g+D_gY)^*(C_g+D_gY)-Y^*Y
=(I+C_g^{-1}D_gY)^*C_g^*C_g(I+C_g^{-1}D_gY)-Y^*Y, 
\]
where $M^*= {^t}{\bar M}$. 

For Hermitian matrices $A, B$, we write
$A\geq B$, if $A-B$ is positive semi-definite, and write 
$A> B$, if $A-B$ is positive definite. 

Note that $G\geq 0$ by (\ref{YCD}), and that 
$\det(I+C_g^{-1}D_gY)\neq 0$
holds for any $Y$ with $\|Y\|<\min\{R^{-1}, R_0^{-1}\}$ 
and any $g\in\Gamma^*$ by Lemma \ref{BdABCD}(iii).
Therefore, the inequality
\[ 
C_g^*C_g \geq {(I+C_g^{-1}D_gY)^*}^{-1}Y^*Y(I+C_g^{-1}D_gY)^{-1}
\]
holds for $\|Y\|<\min\{R^{-1}, R_0^{-1}\}$ and $g\in\Gamma^*$.
Set $Y=tI$, $t=\frac{1}{2}\min\{1, R^{-1}, R_0^{-1}\}$. 
Then, we have
\begin{equation}
 C_g^*C_g \geq 
t^2{(I+tC_g^{-1}D_g)^*}^{-1}(I+tC_g^{-1}D_g)^{-1} 
> \frac{t^2}{4} I.
\end{equation}
Thus we obtain
\[ \|C_g^{-1}\|\leq \frac{2}{t} . \]
\ \qed

\begin{lem}\label{SLIM} 
Suppose that $\Gamma$ satisfies $(\clubsuit)$. 
Then, for any normal sequence $\{g_n\}\subset \Gamma$,  
the sequence $\{\Sigma_{g_n}\}$ 
converges as sets to a single limit $n$-plane,  if and only if  
$\Gamma$ satisfies $(\spadesuit)$. 
\end{lem}
{\bf Proof} \  
Set $\tg_n=\MAT{A_n}{B_n}{C_n}{D_n}\in\SL_{2n+2}(\C)$. 
By the defining equation $\|C_nz'+D_nz''\|=\|z''\|$, 
we have
\[
\Sigma_{g_n}=\{[z':z'']\in\P^{2n+1}\ :\ z'=(C^{-1}_nU-C^{-1}_nD_n)\eta, \ 
z''=\eta, \ \ \eta\in S^{2n+1}, \ U\in \Uni_{n+1} \},
\]
where $S^{2n+1}$  is the unit sphere in $\C^{n+1}$, and 
$\Uni_{n+1}$ is the group of unitary matrices of size $n+1$.
The sequence $\{C_n^{-1}D_n\}$ is 
bounded by Lemma \ref{BdABCD}(iii), and so is 
$\{C_n^{-1}\}$ by the assumption $(\clubsuit)$.
Now consider any subsequence of $\{g_n\}$ such that
$\{C_n^{-1}\}$ converges. Then take again a subsequence such that
$\{C_n^{-1}D_n\}$ also converges. 
Put $L=-\lim_{n\to\infty}C_n^{-1}D_n$ and 
$G=\lim_{n\to\infty} C_n^{-1}$. Then we see
that the set of points on $\Sigma_{g_n}$ converges to 
the set of points
\[
\Sigma := \{ [z':z''] \in\P^{2n+1} : z'=(GU+L)\eta, \ z''=\eta, \ \ \eta\in S^{2n+1}, \ U\in \Uni_{n+1} \}. 
\]
Thus $\Sigma$ consists of a single $n$-plane if and only if $G=0$. 
Here $z'=Lz''$ is the limit $n$-plane of $\{g_n^{-1}(\ell'')\}$.  
This implies the lemma. 
\ \qed

\begin{lem}\label{aInSigma}
Suppose that $\Gamma$ satisfies $(\clubsuit)$ and $(\spadesuit)$.  
Then, for $a\in \Omega$, there can 
be at most finite number of $\Sigma_g$
which contains $a$. 
\end{lem}
{\bf Proof} \  
Suppose that there is an infinite number of $g_n\in\Gamma$, 
$n=1,2,\dots$, such that $a\in\Sigma_{g_n}$.
Then taking a normal subsequence of $\{g_n\}$, we see 
that $a$ is on a limit $n$-plane by Lemma \ref{SLIM}, 
since $\Gamma$ satisfies $(\spadesuit)$.
This contradicts $a\in\Omega$. 
\ \qed

\bigskip
Now recall the definition of $\mu_g$ for $g\in\Gamma^*$. 
We also define 
\[ \mu_1(z) \equiv 1 \ \ \ \mbox{for} \ \ \ g=1. \]

\begin{lem}\label{mu_eq} For $g, h \in \Gamma$, we have 
\[ 
\mu_{h\circ g}(z) = \mu_h(g(z))\mu_g(z), 
\ \ \ z \in \P^{2n+1} \setminus \{g^{-1}(\ell'')\cup (h\circ g)^{-1}(\ell'')\}. 
\]
\end{lem}
{\bf Proof} \  Easy by Lemma \ref{mu}. 
\ \qed

\begin{lem} \label{Ggc} For any $g \in \Gamma$, 
we have $g(\bar\Delta) \subset \Delta_{g^{-1}}^c$. 
\end{lem}
{\bf Proof} \  By Lemma \ref{mu_eq}, we have 
$\mu_{g^{-1}}(g(z))=\mu_g(z)^{-1}$. 
Since $\mu_g(z)\leq 1$ for $z \in \bar\Delta$, we have
the lemma. 
\ \qed

\begin{thm} \label{Ford_Rgn} 
Let $\Gamma\subset \PGL_{2n+2}(\C)$ 
be a type $\ELL$ group. Assume that $\Gamma$ is torsion free and 
satisfies both $(\clubsuit)$ and $(\spadesuit)$.
Then $F$ has the following properties :
\begin{enumerate}
\item[(1) ] For $g \in \Gamma$, $g(F)\subset F$ holds if and only if $g=1$.
\item[(2) ] For $g \in \Gamma^*$, 
$g(F)\cap F=\emptyset$ holds.
\item[(3) ] For every $z \in \Omega$, 
there is an element $g \in \Gamma$
such that $g(z) \in \Omega\cap\bar F$.
\item[(4) ] Suppose that the equality $w=g(z)$  holds
for some $z, w \in \Omega\cap\bar F$ and $g \in \Gamma^*$,
then both $z$ and $w$ are on $\Omega\cap\pd\bar F$.
\end{enumerate}
\end{thm}
{\bf Proof} \  
The following proof is an analogue of \cite[pp.33-34]{Maskit}. 

$(1)$ \ By Proposition \ref{NbdOfL2}, $F$ contains a 
tubular neighborhood $W$ of $\ell''$.
Suppose that $g(F)\subset F$, Then, by Lemma \ref{mu_eq}, 
we have 
\[ 
\mu_{h\circ g}(z)=\mu_h(g(z))\mu_g(z) \leq \mu_g(z) \leq 1
\]
on $W$ for any $h$. 
Letting $h=g^{-1}$, we see that $\mu_g(z)=1$ 
on $W$. This implies $C_g=0$. 
Hence $g(\ell'')=\ell''$. Since $\ell''$ is not a limit $n$-plane, we see that
$g$ is of finite order. Since $\Gamma$ is torsion free by assumption,
we see that $g=1$. The converse is obvious.

$(2)$ \ By Lemma \ref{Ggc}, we have 
$g(\bar\Delta) \subset \Delta_{g^{-1}}^c$. 
This implies $g(\bar\Delta)\cap \Delta_{g^{-1}} = \emptyset$.
Hence $g(F)\cap F\subset g(F)\cap\bar\Delta\subset g(F)\cap\bar\Delta_{g^{-1}}\subset g(\bar\Delta)\cap \Delta_{g^{-1}}=\emptyset$.

$(3)$ \ Take a point $z$ in $\Omega$. 
If $g(z) \in \ell''$ for some $g \in\Gamma$, 
then $g(z) \in \Omega\cap \bF$ by the assumption $(\clubsuit)$ and 
Proposition \ref{NbdOfL2}. 
Therefore we can assume that $g(z)\notin\ell''$ for any $g$. 
Then $\mu_g(z)$ is defined and hence has a finite-value for any $g$. 
By the assumptions $(\clubsuit)$ and $(\spadesuit)$, 
 $\mu_g(z)<1$ holds for all 
except for finitely many $g \in \Gamma$. 
Therefore we can choose $g$ such that $\mu_g(z)$ is maximal among
all $g$. 
Then, by Lemma \ref{mu_eq}, 
we have $\mu_h(g(z))\leq 1$ for any $h \in \Gamma$.
This implies $g(z) \in\bar\Delta$.
Thus we have $\Omega\subset \bigcup_{g\in\Gamma}g(\bar\Delta)$
and hence 
\begin{equation}\label{ogd}
\Omega = \bigcup_{g\in\Gamma}g(\Omega\cap\bar\Delta).
\end{equation}
We claim that the set $\Omega\cap(\bar\Delta\setminus\bar F)$ is empty.
To verify this, we suppose contrary that 
a point $w\in\Omega\cap(\bar\Delta\setminus\bar F)$ exists.
Since $\bar\Delta\setminus\bar F$ is thin in $\P^{2n+1}$, so is  
$\bigcup_{g\in\Gamma}g(\Omega\cap(\bar\Delta\setminus \bar F))$.
Hence, by 
\[
\Omega\setminus
\bigcup_{g\in\Gamma}g(\bar F)
=
\bigcup_{g\in\Gamma}g(\Omega\cap\bar\Delta)\setminus
\bigcup_{g\in\Gamma}g(\Omega\cap\bar F)
\subset 
\bigcup_{g\in\Gamma}g(\Omega\cap(\bar\Delta\setminus\bar F)),
\]
we see that the set $\Omega\setminus\bigcup_{g\in\Gamma}g(\bar F)$
is thin in $\Omega$. Therefore, we can find sequences 
$\{w_n\}\subset \Omega\cap \bar F\subset \bDelta$
and $\{g_n\}\subset\Gamma$ 
such that $\lim_{n\to\infty}g_n(w_n)=w$.
Since $w\notin\pd\bar F$, $\{g_n\}$ can be chosen to be 
a sequence of distinct elements. 
By Lemma \ref{BdABCD} and the assumptions $(\clubsuit)$ and $(\spadesuit)$, 
we can choose a subsequence 
of $\{g_n\}$ such that 
the $n$-plane $g_n(\ell'')=\{z'+C_{g_n^{-1}}^{-1}D_{g_n^{-1}}z''=0\}$
converges to a limit $n$-plane 
$\ell_L=\{z'=Lz''\}$, 
$L=-\lim_{n\to\infty}C_{g^{-1}_n}^{-1}D_{g^{-1}_n}$, 
and such that $\lim_{n\to\infty}\|C_{g_n^{-1}}^{-1}\|=0$ holds.
The set $\Delta_{g_n^{-1}}^c$ is the image of the map
\[ \vf_n : [0,1]\times U_{n+1}\times S^{2n+1} \to \P^{2n+1} \]
defined by 
\[ 
\vf_n : (t, U, \eta) \mapsto 
\left[t C_{g_n^{-1}}^{-1}U\eta - C_{g_n^{-1}}^{-1}D_{g_n^{-1}}\eta : \eta\right].
\]
Therefore, the sequence $\{\Delta_{g_n^{-1}}^c\}$ of sets converges to the 
image of the limit map
\[ \vf : [0,1]\times U_{n+1}\times S^{2n+1} \to \P^{2n+1},
 \ \ \ (t, U, \eta) \mapsto [L\eta : \eta],
\]
which is $\ell_L$. 
Since $g_n(\ell'')\subset g_n(\bar\Delta)\subset \Delta_{g_n^{-1}}^c$
holds by Lemma \ref{Ggc},
$\{g_n(\bDelta)\}$ also converges to $\ell_L$. 
Hence $w$ is on the limit $n$-plane $\ell_L$.
Since $w \in \Omega$,  this is absurd. 
Thus our claim is verified. 
Now, by $(\ref{ogd})$, we have
\begin{equation}\label{ogF}\nonumber
\Omega = \bigcup_{g\in\Gamma}g(\Omega\cap\bar F).
\end{equation} 
This proves $(3)$.

$(4)$ \ By $(2)$, either $z\in\pd\bar F$ or $w\in\pd\bar F$. 
Replacing $g$ with $g^{-1}$ if necessary, 
we can assume that $z\in\pd\bar F$. 
Since $z, w \in\bar F$, $\mu_f(z)\leq 1$ and $\mu_f(w)\leq 1$
hold for any $f\in\Gamma^*$. 
Hence, by the equality 
\[
1=\mu_1(z)=\mu_{g^{-1}g}(z)=\mu_{g^{-1}}(g(z))\mu_g(z)
=\mu_{g^{-1}}(w)\mu_g(z),
\]
we have $\mu_{g^{-1}}(w)=\mu_g(z)=1$. Hence, in particular,  
we have $w \in \Sigma_{g^{-1}}$. 
On the other hand, there can be at most finite number of 
$\Sigma_f$ with $w\in\Sigma_f$ by Lemma \ref{aInSigma}.
This implies $w\in\pd\bar F$. 
\ \qed

\begin{rem} For type $\ELL$ groups, both conditions $(\clubsuit)$ and 
$(\spadesuit)$ are automatically satisfied, if the series
\[ \sum_{g\in\Gamma^*}\|C_g^{-1}\|^\delta \]
is convergent for some constant $\delta>0$. 
\end{rem}
Theorem \ref{Ford_Rgn} is useful when we check 
the quotient space $\Omega(\Gamma)/\Gamma$ becomes compact or not.
See examples in subsection \ref{three}.

\section{Examples}
\subsection{Type $\ELL$ groups in general dimension}\label{general}
Suppose that we are given two groups $\Gamma_\nu\subset\PGL_{2n+1}(\C)$, 
$\nu=1,2$, of type $\ELL$ which satisfy $(\clubsuit)$ and $(\spadesuit)$.
Applying a Klein combination, we can construct another group 
$\Gamma$ of type $\ELL$ which is isomorphic to the free product 
$\Gamma_1*\Gamma_2$. 
In this subsection, we show that, 
replacing $\Gamma_\nu$ with their suitable conjugate subgroups in 
$\PGL_{2n+1}(\C)$, we can make $\Gamma$ also satisfy 
both $(\clubsuit)$ and $(\spadesuit)$.

\bigskip
Let $F_\nu$ the $F$-region of $\Gamma_\nu$ with respect to $[z]=[z':z'']$.
Let $\rho_\nu>0$ be the numbers such that $\|C_{g}^{-1}\|\leq\rho_\nu$
for all $g \in \Gamma^*_\nu$.

\begin{lem}\label{ra} 
Let $a\in\R$ be a positive constant, and consider the new system of coordinates 
$[\zeta':\zeta'']$ on $\P^{2n+1}$ defined by
\[ \zeta'=az', \ \ \zeta''=a^{-1}z''. \]
Let $\al=\MAT{aI}{0}{0}{a^{-1}I}\in\PGL_{2n+2}(\C)$. 
Then, for a given number $r>0$, 
we can choose $a>0$ so that the following are satisfied simultaneously.
\begin{enumerate}
\item[(i) ] $\al(F_1\cap F_2)$ contains 
the set $V=\{\|\zeta'\|\geq r\|\zeta''\|\}$.
\item[(ii) ] $\|C_{g}^{-1}\|\leq 1$
for all $g=\MAT{A_g}{B_g}{C_g}{D_g} \in \al\Gamma^*_\nu\al^{-1}$, 
$\nu=1,2$.
\end{enumerate}
\end{lem}
{\bf Proof} \  
Since $F_\nu$ contains a tubular neighborhood of $z''=0$, there 
is $r_1>0$ such that
\[ \{\|z'\|\geq r_1\|z''\|\}\subset F_1\cap F_2. \] 
Choose $a>0$ satisfying 
\begin{equation}\label{r1}
 a^2 \leq r_1^{-1}r.
\end{equation}
Take any $[\zeta':\zeta'']\in V$, and set $[z':z'']=\al^{-1}([\zeta':\zeta''])$.
Then $z'=a^{-1}\zeta'$ and $z''=a\zeta''$, and we have 
\[
 \|z'\|=a^{-1}\|\zeta'\|\geq a^{-1}r\|\zeta''\|=a^{-2}r\|z''\|\geq r_1\|z''\|.
\]
Hence $[z':z'']\in F_1\cap F_2$. 
This shows that $[\zeta':\zeta'']\in\al(F_1\cap F_2)$. 
Thus  (i) is satisfied for $a>0$ with (\ref{r1}).

Let $\gm=\MAT{A}{B}{C}{D}\in\Gamma_\nu^*$. Then by
$g=\al\gm\al^{-1}=\MAT{A}{a^2B}{a^{-2}C}{D}$, 
we have $\|C_g^{-1}\|=a^2\|C^{-1}\|$. Therefore the number 
$a>0$ with 
\begin{equation}\label{r2}
a^2\leq\rho_\nu^{-1}, \ \ \nu=1,2
\end{equation} 
satisfies (ii).  
Thus it is enough to choose $a>0$ which satisfies (\ref{r1}) and (\ref{r2}). 
\ \qed
 
\bigskip
Fix $a>0$ such that (i) and (ii) in the lemma above hold
and replace the original coordinates $[z':z'']$ 
with $[\zeta':\zeta'']$, and $\Gamma_\nu$ with $\al\Gamma_\nu\al^{-1}$,
and $F_\nu$ with $\al(F_\nu)$. We use the original notation such as
$[z':z'']$, $\Gamma_\nu$, and $F_\nu$ to avoid abuse of notation.
Let $U_r=\{\|z'\|\leq r\|z''\|\}$. Then $\gm(F_\nu)\subset U_r$ for 
$\gm\in\Gamma^*_\nu$, $(\nu=1,2)$.

Put $\sigma =\MAT{0}{I}{I}{0} \in \PGL_{2n+2}(\C)$,
and consider the sets $F_1$ and $\sigma(F_2)$ as two subsets in the 
same projective space $\P^{2n+1}$. 
Note that the set $\{r\|z''\|\leq \|z'\|\leq r^{-1}\|z''\|\}$, 
$0<r<1$, is contained in $F_1\cap\sigma(F_2)$, and that 
$\gm(\sigma(F_2))\subset \sigma(U_r)$ for 
$\gm\in\sigma\Gamma_2\sigma^{-1}$.

Put 
$\tau =\MAT{I}{I}{-I}{I}\in\PGL_{2n+2}(\C)$.
Introduce a new coordinate system $[w]=[w':w'']$ by $w=\tau(z)$. 
Put $\Gamma_1'=\tau\Gamma_1\tau^{-1}$ and 
$\Gamma_2'=\tau(\sigma\Gamma_2\sigma^{-1})\tau^{-1}$.
Let $\Gamma$ be the group generated by $\Gamma_1'$ and $\Gamma_2'$.
This subsection is devoted to prove the following theorem. 
\begin{thm}\label{KC}
$\Gamma$ is a group of type $\ELL$ which satisfy 
$(\clubsuit)$ and $(\spadesuit)$ with respect to $[w]$. 
\end{thm}
{\bf Proof} \  
By the construction, we see that $\tau(F_1\cap\sigma(F_2))$
is a fundamental set of $\Gamma$.
The $n$-plane $\{w''=0\}$ has a tubular neighborhood contained in 
$\tau(F_1\cap\sigma(F_2))$. Thus, by Proposition \ref{NbdOfL2}, 
it is enough to show that $(\spadesuit)$ is satisfied.
 
\begin{lem}\label{PQ}
Let $\ell_P : w'=Pw''$ and $\ell_Q : w'=Qw''$ be 
$n$-planes in $\tau(U_r)$ and $\tau\sigma(U_r)$,
respectively. Then there is a positive constant $K_r$ such that 
\[ \|(P-Q)^{-1}\|\leq K_r, \]
where 
\[
\lim_{r\to 0}K_r = 2^{-1}.
\]
\end{lem}
{\bf Proof} \  
Take the $n$-planes 
$\ell_X : z'=Xz''$ in $U_r$ and 
$\ell_Y : z''=Yz'$ in $\sigma(U_r)$ 
such that $\tau(\ell_X)=\ell_P$ and $\tau(\ell_Y)=\ell_Q$, respectively.
Then we have
\[ P=(I+X)(I-X)^{-1}, \ \ \ Q=-(I+Y)(I-Y)^{-1}. \]
Since $\ell_P \cap \ell_Q=\emptyset$, $\det(P-Q)\neq 0$ holds.
Hence
\[ 
\|(P-Q)^{-1}\|
=\left\|\left((I+X)(I-X)^{-1}+(I+Y)(I-Y)^{-1}\right)^{-1}\right\|. 
\]
Set
\[ 
K_r=\sup_{\{\|X\|\leq r, \ \|Y\|\leq r\}}
\left\|\left((I+X)(I-X)^{-1}+(I+Y)(I-Y)^{-1}\right)^{-1}\right\|
\]
It is clear that $K_r$ is finite for $0\leq r<1$, 
and that $\lim_{r\to 0}K_r=2^{-1}$. 
This implies the lemma. 
\ \qed

\bigskip
Any element $f\in\Gamma^*$ can be written in the following {\em normal}
 form \footnote{\cite[p.136]{Maskit}}
\[ f=g_m\cdots g_1. \]
Here either $g_{2j+1}\in{\Gamma'}^*_1$, $g_{2j}\in{\Gamma'}^*_2$, 
or $g_{2j+1}\in{\Gamma'}^*_2$, $g_{2j} \in{\Gamma'}^*_1$.
The number $m$ is called the {\em length} of $f$, which is denoted by $|f|$.

\begin{lem}\label{Cggg} Take any element $f\in\Gamma^*$, and  
write $f$ in the normal form :
\[ f=g_m\cdot g_{m-1}\cdots g_1. \]
Then we have 
\begin{equation}\label{CCC}
 \|C_f^{-1}\|\leq 
K_r^{m-1}\prod_{j=1}^m\|C_{g_j}^{-1}\|, 
\end{equation}
where 
\[
f=\MAT{A_f}{B_f}{C_f}{D_f}, \ 
g_j=\MAT{A_{g_j}}{B_{g_j}}{C_{g_j}}{D_{g_j}}.
\]
\end{lem}
{\bf Proof} \   
Set $g=g_m$ and $h=g_{m-1}\cdots g_1$. 
Comparing the components of $f=gh$, we have
\begin{equation}\label{fgh}
 C_f=C_gA_h+D_gC_h=C_g(A_hC_h^{-1}+C_g^{-1}D_g)C_h. 
\end{equation}

First assume that $g_1\in{\Gamma'}^*_1$. 
If $g\in{\Gamma'}^*_1$, then $|f|=m$ is odd. 
Since $g^{-1}\in\Gamma^*_1$, 
$g^{-1}(\{w''=0\})=\{w'=-C_g^{-1}D_gw''\}\subset\tau(U_r)$. 
Since $|h|$ is even, we see that 
$h(\{w''=0\})=\{w'=A_hC_h^{-1}w''\} \subset\tau\sigma(U_r)$.
Since $\tau(U_r)\cap\tau\sigma(U_r)=\emptyset$,  we see that 
$\det(A_hC_h^{-1}+C_g^{-1}D_g)\neq 0$. 
Hence $C_f$ is also non-singular, and by (\ref{fgh}), we have
\begin{equation}\label{Cghinv}
C_f^{-1}=C_h^{-1}(A_hC_h^{-1}+C_g^{-1}D_g)^{-1}C_g^{-1}.
\end{equation}
By Lemma \ref{PQ}, it follows that
\begin{equation}\label{PQAC}
\|(A_hC_h^{-1}+C_g^{-1}D_g)^{-1}\|\leq K_r.
\end{equation}
Hence by (\ref{Cghinv}), we have
\begin{equation}\label{KCC}
\|C_f^{-1}\| \leq K_r \|C_g^{-1}\|\cdot\|C_h^{-1}\|.
\end{equation}

If $g\in{\Gamma'}^*_2$, then $|f|=m$ is even. 
Since $g^{-1}\in\Gamma^*_2$, 
$g^{-1}(\{w''=0\})=\{w'=-C_g^{-1}D_gw''\}\subset\tau\sigma(U_r)$. 
Since $|h|$ is odd, we see that 
$h(\{w''=0\})=\{w'=A_hC_h^{-1}w''\} \subset\tau(U_r)$.
Then by the same argument as the case $g=g_m\in{\Gamma'}^*_1$,
we obtain (\ref{KCC}).

\bigskip
Next assume that $g_1\in{\Gamma'}^*_2$. 
If $g\in{\Gamma'}^*_1$, then $|f|=m$ is even. 
Since $g^{-1}\in\Gamma^*_1$, 
$g^{-1}(\{w''=0\})=\{w'=-C_g^{-1}D_gw''\}\subset\tau(U_r)$. 
Since $|h|$ is odd, we see that 
$h(\{w''=0\})=\{w'=A_hC_h^{-1}w''\} \subset\tau\sigma(U_r)$.
Then the rest of the argument is the same as above, and we obtain (\ref{KCC}).

If $g\in{\Gamma'}^*_2$, then $|f|=m$ is odd. 
Since $g^{-1}\in\Gamma^*_2$, 
$g^{-1}(\{w''=0\})=\{w'=-C_g^{-1}D_gw''\}\subset\tau\sigma(U_r)$. 
Since $|h|$ is even, we see that 
$h(\{w''=0\})=\{w'=A_hC_h^{-1}w''\} \subset\tau(U_r)$.
Then the rest of the argument is the same as above, and we obtain (\ref{KCC}).

The lemma follows from (\ref{KCC}) by induction on $m$. 
\ \qed

\bigskip\noindent
{\bf Proof of Theorem \ref{KC} (continued)}
It remains to show that $\Gamma$ satisfies $(\spadesuit)$. 
By Lemma \ref{ra}, we can assume that $\rho_\nu\leq 1$, $\nu=1,2$.
By Lemma \ref{PQ}, we fix small $r$, $0<r<1$, such that $K_r<1$ holds.
Now we shall show that $\Gamma$ satisfies $(\spadesuit)$.

Suppose that $(\spadesuit)$ does not hold. Then there is a sequence 
$\{f_m\}_m\subset\Gamma^*$ such that 
\begin{equation}\label{ep}
 \lim_{m\to\infty}\|C_{f_m}^{-1}\|=\ep > 0.
\end{equation}
If there is a subsequence $\{h_m\}$ of $\{f_m\}$ such that
$\lim_{m\to\infty}|h_m|=\infty$. Then, 
$\lim_{m\to\infty}\|C_{h_m}^{-1}\|=0$ follows 
from $K_r<1$ and $\rho_\nu\leq 1$ by Lemma \ref{Cggg}.
This contradicts (\ref{ep}). Therefore the sequence $\{|f_m|\}_m$ of lengths
is bounded.  Let $b$ be a bound of $\{|f_m|\}_m$, i.e, 
$|f_m|\leq b$ for all $m$. Write $f_m$ in the "extended" normal form, 
\[ f_m=g_{m,b}g_{m,b-1}\cdots g_{m,1},\]
where $g_{m,|f_m|}\cdots g_{m,1}$ is the normal form of $f_m$, and 
$g_{m,j}=1$ for $|f_m|<j\leq b$. 
Since both $\Gamma^*_1$ and $\Gamma^*_2$ satisfy
$(\clubsuit)$ and $(\spadesuit)$, 
we can find some $k$, $1\leq k\leq b$, such that 
$\{\|C_{g_{m,k}}^{-1}\|\}_m$ contains a subsequence which 
converges to zero. This implies that the corresponding subsequence of 
$\{\|C_{f_m}^{-1}\|\}_m$ also converges to zero. 
This again contradicts (\ref{ep}). 
\ \qed

\begin{rem} A typical higher dimensional example treated in this subsection 
is a Schottky group. Let $\Gamma$ be the infinite cyclic group generated by
$g=\MAT{A}{0}{0}{B}\in\PGL_{2n+2}(\C)$. Let $\al_j$ be
the eigen-values of $A$, and $\be_k$ the eigen-values of $B$.
Assume that $|\al_j| < |\be_k|$ holds for any pairs $(\al_j, \be_k)$.
Then $\Gamma$ is a type $\ELL$ group, 
where $\Omega(\Gamma)=\P^{2n+1}\setminus(\{z'=0\}\cup\{z''=0\})$.
Introduce a new coordinate $[w':w'']$ by $w'=z'+z''$ and $w''=-z'+z''$.
Then $\Gamma$ satisfies $(\clubsuit)$ and $(\spadesuit)$ with respect to
$[w':w'']$. By successive Klein combinations, 
we can get type $\ELL$ groups with $(\clubsuit)$ and $(\spadesuit)$ 
with respect to some coordinate system.
\end{rem}

\subsection{Type $\ELL$ groups in dimension $3$}\label{three}
In this subsection, we shall give three examples of type $\ELL$ groups.
If a finitely generated discrete infinite subgroup $\Gamma \subset \PGL_4(\C)$ 
admits an invariant surface $S$ in $\P^3$ and never admits invariant planes, 
then $S$ is necessarily one of the following: 
(i) the tangential surface of a twisted cubic curve, 
(ii) a non-singular quadric surface, 
(iii) a cone over a non-singular conic \cite{plane}. 
Each case has examples of type $\ELL$ groups with 
$(\clubsuit)$ and $(\spadesuit)$. 
In the cases (i) and (ii), there are examples with compact connected 
canonical quotients. 
The example for the case (ii) is due to Fujiki\cite{Hayama}. 
For the case (iii), we have only one example at present, 
whose canonical quotient is connected and non-compact,
but it has an invariant plane.

\subsubsection{Kleinian groups acting on a twisted cubic curve}\label{tw}
Fix a twisted cubic curve $C \subset \P^3$, 
which is defined to be the image of the map
\[ \tau : \P^1 \to \P^3, \ \ \ \tau([s:1]) = [s^3:s^2:s:1]. \]
Then $\tau$ determines a group representation
\[ \tau_* : \PSL_2(\C) \to \PSL_4(\C) \]
such that $\tau\circ g = \tau_*(g)\circ\tau$.
Explicitly, for 
$ g=\pm\MAT{a}{b}{c}{d} \in \PSL_2(\C)$, 
$\tau_*(g)$ is given by
$ \tau_*(g) = \pm \ttau_*(g) \in \PSL_4(\C)$, 
where 
\begin{equation}\label{tgamma}
\ttau_*(g) = \left(
\begin{array}{cccc}
a^3 & 3a^2b & 3ab^2 & b^3 \\ 
a^2c & a^2d+2abc & 2abd+b^2c & b^2d \\
ac^2 & 2acd+bc^2 & ad^2+2bcd & bd^2 \\
c^3 & 3c^2d & 3cd^2 & d^3
\end{array}
\right) \in \SL_4(\C).   
\end{equation}

The group $\PGL_4(\C)$ acts on the set of lines in $\P^3$, i.e., on 
$\Grass$. Using Pl\"ucker coordinates, we can embed $\Grass$ into 
$\P^5=\P(\wg^2\C^4)$. Since any $A\in \GL_4(\C)$ defines a linear automorphism
on $\wg^2\C^4\simeq\C^6$, we have the group homomorphism
\[ \rho : \GL_4(\C) \to \GL_6(\C). \]
Let $e_0={^t(1,0,0,0)},  e_1={^t(0,1,0,0)},  e_2={^t(0,0,1,0)},  e_3={^t(0,0,0,1)}$,
and $e_j \wg e_k$ the linear 2-space spanned by $\{e_j,e_k\}$,
where $e_j\wedge e_k = - e_k\wedge e_j$. 
In this subsection in the following, we write 
$\tg = \rho\circ\ttau_*(g) \in \GL_6(\C)$, 
which is well-defined for $g=\pm\MAT{a}{b}{c}{d}\in\PSL_2(\C)$. 
Then, with respect to the basis 
\[ 
\{ e_0\wg e_1, \ e_0\wg e_2, \ e_0\wg e_3, 
  \ e_1\wg e_2, \ e_1\wg e_3, \ e_2\wg e_3\},
\]
in $\C^6$, $\tg=\rho\circ\ttau_*(g) \in \GL_6(\C)$ is given by
\begin{small}
\begin{equation}\label{size6}\nonumber
\tg=\left(
\begin{array}{cccccc}
a^4 & 2a^3b & a^2b^2 &  3a^2b^2 & 2ab^3 & b^4 \\
2a^3c & a^2(ad+3bc) & ab(ad+bc) & 3ab(ad+bc) & b^2(3ad+bc) & 2b^3d\\
3a^2c^2 & 3ac(ad+bc) & a^2d^2+abcd+b^2c^2 & 9abcd & 3bd(ad+bc) & 3b^2d^2\\
a^2c^2 & ac(ad+bc) & abcd &a^2d^2 + abcd+ b^2c^2 & bd(ad+bc) &b^2d^2\\
2ac^3 & c^2(3ad+bc) & cd(ad+bc) & 3cd(ad+bc) & d^2(ad+3bc) &2bd^3\\
c^4 & 2c^3d & c^2d^2 & 3c^2d^2 & 2cd^3 &d^4
\end{array}
\right)
\end{equation}
\end{small}

\bigskip
\underline{\bf Limit sets} 
In the following in this subsection, we 
let $\Gamma \subset \PSL_2(\C)$ be a Kleinian group whose 
set of discontinuity $\Omega_{\P^1}$ contains $[1:0]\in\P^1$.
Put $\Lambda_{\P^1}=\P^1\setminus\Omega_{\P^1}$.
We consider the group 
$\tGamma=\tau_*(\Gamma)$, 
which we regard as a subgroup of $\PGL_4(\C)$.

\bigskip
Let $\{\gamma_n\} \subset \Gamma$ be a normal sequence. 
Let  
\[ 
g_n = \MAT{a_n}{b_n}{c_n}{d_n} \in \SL_2(\C), 
\ \ \ n=1,2,\dots
\]
be a sequence of representatives of $\{\gamma_n\}$ 
such that $\{c_n^{-1}g_n\}$
converges to a matrix of the form
$ h=\MAT{\mu}{-\lm\mu}{1}{-\lm} \in M_2(\C)$, $\lm, \mu\in\C$, 
since $\{a_nc_n^{-1}\}$ and $\{c_n^{-1}d_n\}$ are 
bounded (cf. Lemma \ref{BdABCD}). Put 
\[ G_n=c_n^{-4}\ttau_*(g_n) \in \GL_6(\C).\] 
Then
\begin{small}
\begin{equation}\label{rhogamma}\nonumber
G:=\lim_n G_n=
\left(
\begin{array}{cccccc}
\mu^4  & -2\lm\mu^4 & \lm^2\mu^4 &  3\lm^2\mu^4 & -2\lm^3\mu^4 & \lm^4\mu^4 \\
2\mu^3 & -4\lm\mu^3 & 2\lm^2\mu^3 & 6\lm^2\mu^3 & -4\lm^3\mu^3 & 2\lm^4\mu^3 \\
3\mu^2 & -6\lm\mu^2 & 3\lm^2\mu^2 & 9\lm^2\mu^2 & -6\lm^3\mu^2 & 3\lm^4\mu^2 \\
\mu^2  & -2\lm\mu^2 & \lm^2\mu^2 & 3\lm^2\mu^2 & -2\lm^3\mu^2 & \lm^4\mu^2 \\
2\mu & -4\lm\mu & 2\lm^2\mu & 6\lm^2\mu & -4\lm^3\mu & 2\lm^4\mu \\
1 & -2\lm & \lm^2 & 3\lm^2 & -2\lm^3 & \lm^4
\end{array}
\right). 
\end{equation}
\end{small}
The limit $G$ defines a projection to the limit image
\[ \P^5\setminus H \to I:=\{ [\mu^4:2\mu^3:3\mu^2:\mu^2:2\mu :1] \}. \]
The limit kernel $H$ is the 4-plane defined by 
\[ \{ \zeta=[\zeta_j] \in \P^5 : 
\zeta_0 -2\lm\zeta_1+\lm^2\zeta_2 + 3\lm^2\zeta_3 
- 2\lm^3\zeta_4 + \lm^4\zeta_5 =0 \}. 
\]
Let $\ell_\mu$ the tangent line to the curve $C$ 
at $[\mu:1]$. Then $\hat\ell_\mu \in \Grass \subset \P^5$ 
is given by
\[
\left(3\mu^2e_0 + 2\mu e_1+ e_2\right)
\wg\left(\mu^3e_0 + \mu^2e_1 + \mu e_2 + e_3\right) 
\]
\[
= \mu^4 e_0\wg e_1 + 2\mu^3 e_0\wg e_2 +3\mu^2 e_0\wg e_3 +\mu^2 e_1\wg e_2 
+ 2\mu e_1\wg e_3 + e_2\wg e_3
\]
which is nothing but 
the limit image $I=[\mu^4:2\mu^3:3\mu^2:\mu^2:2\mu:1] \in \P^5$. 
Hence $\ell_\mu$ is the limit image of the sequence $\{\tau_*(\gamma_n)\}$. 
Here the limit kernel $H\cap\Grass$ is the set of lines in $\P^3$ which intersect 
the tangent line to $C$ at the limit point $\tau([\lm : 1])$.  
Thus we have the following result.

\begin{thm}\label{limset}  
Let $\Gamma \subset \PSL_2(\C)$ be a Kleinian group. 
Then 
\[ \tGamma=\tau_*(\Gamma) \subset \PGL_4(\C) \]
is a group of type $\ELL$. The limit set is given by 
\[ \Lambda(\tGamma) 
= \bigcup_{\lm \in \Lambda_{\P^1}}\left|\ell_\lm\right|,
\]
where 
$\left|\ell_{\lm}\right|$ is the support of the tangent line $\ell_\lm$ 
to the twisted cubic curve at $\tau([\lm:1])$. 
\end{thm}

\begin{prop}\label{P_tw} \nonumber
Let $\Gamma \subset \PSL_2(\C)$ be a Kleinian group whose set of
discontinuity contains the point $[1:0]\in \P^1$.
Then the series
\begin{equation}\label{PoincareTw} 
 \mathop{\sum_{\tg \in \tGamma^*}}
\|C_g^{-1}\|^\delta 
\end{equation}
is convergent for any $\delta\geq 4$. Thus $\Gamma$ satisfies 
$(\clubsuit)$ and $(\spadesuit)$. 
\end{prop}
{\bf Proof} \  
By our choice of coordinates on $\P^1$, we see that
$c_g\neq 0$ for $g \neq 1$, and $\{a_g/c_g\}_g, \{b_g/c_g\}_g, \{d_g/c_g\}_g$
are uniformly bounded. 
Since
\[  
C_g =\left(
\begin{array}{cc}
a_gc_g^2 & 2a_gc_gd_g+b_gc_g^2 \\
c_g^3 & 3c_g^2d_g 
\end{array}\right)  
=
c_g^3\left(\begin{array}{cc}
a_g/c_g & 2a_gd_g/c_g^2 + b_g/c_g \\
1 & 3d_g/c_g 
\end{array}\right),
\]
we see that
\begin{equation}\label{Cgcg}\nonumber
 \|C_g^{-1}\| \leq M|\det C_g|^{-1}|c_g|^3 = M|c_g|^{-1} 
\end{equation}
holds  for some $M>0$. 
It is well-known that, for Kleinian groups, the series
$\sum_{g \in \Gamma^*} |c_g|^{-4}$ is convergent
 (\cite[Theorem II.B.5]{Maskit}). 
Hence we have the proposition. 
\ \qed

\bigskip
\underline{\bf Compact quotients}
As an application of Theorem \ref{Ford_Rgn}, we obtain the following.
\begin{thm}\label{TQ} If $\Gamma$ is convex-cocompact
\footnote{geometrically finite and no parabolic elements, 
\cite[p.95]{Kap}} , 
then $\left(\P^3\setminus \Lambda(\tGamma)\right)/\tGamma$ 
is compact. 
\end{thm}
{\bf Proof} \  
Let $\Gamma \subset \SL_2(\C)$ be a convex-cocompact 
Kleinian group. It is known that every limit point of 
a convex-cocompact group is a point of approximation 
\cite[Definitions 4.43, 4.71 and 4.76]{Kap}.
We can assume further that $\Gamma$ is torsion free without loss of 
generality.
By Theorem \ref{Ford_Rgn}, it is enough to show that the set $\bF$
defined by (\ref{truFord})  is 
a compact subset contained in $\Omega(\tGamma)$.  
If $\bDelta$ is contained in $\Omega(\tGamma)$, then 
the quotient $\Omega(\tGamma)/\tGamma$ becomes compact,  
since $\bF \subset \bDelta$ and $\bDelta$ is compact. 
Thus it sufficient to show the following proposition. 
\ \qed

\begin{prop} 
Any limit line does not intersect $\bar{\Delta}$.
\end{prop}
{\bf Proof} \  
Let $\ell_\lm$ be any limit line, 
which is the tangent line to $C$ at $\tau([\lm:1])$, 
$\lm=[\lm:1] \in\Lambda_{\P^1}$.
More explicitly, $\ell_\lm$ is given by
$z' = L_\lm z''$, 
where $z=[z':z''] \in \P^3$ and 
\[ L_\lm = \left(
\begin{array}{cc}
  3\lm^2  & -2\lm^3 \\
 2\lm & -\lm^2 
\end{array} \right).
\]
Recall that every limit point of $\Gamma$ is a point of approximation.
Hence, there are a sequence $\{g_m\}$ of distinct elements 
of $\Gamma$  
and a constant $\delta>0$ such that 
\begin{equation}\label{approx} 
|g_m(\lm) - g_m(\infty)|\geq \delta
\end{equation}
for any $m$. Let 
\[ g_m=
\left(\begin{array}{cc} a_m & b_m \\ c_m & d_m\end{array}\right)
\ \in \SL_2(\C)
\]
The inequality $(\ref{approx})$ is equivalent to
\[
\left|\frac{a_m\lm + b_m}{c_m\lm+d_m} - \frac{a_m}{c_m}\right|
\geq \delta.
\]
This implies
\begin{equation}\label{approx1}
|c_m(c_m\lm+d_m)|\leq \delta^{-1}. 
\end{equation}
Since $\infty=[1:0]\in\Omega_{\P^1}$, 
we know that $\lim_{m\to\infty}|c_m|=\infty$. 
Hence, it follows from (\ref{approx1}) that
\begin{equation}\label{approx2}
\lim_{m\to\infty}|c_m\lm+d_m|=0, \ \ \ 
\lim_{m\to\infty}|c_m(c_m\lm+d_m)^2|=0.
\end{equation}
Again, since $\infty=[1:0]\in\Omega(\Gamma)$, 
there is a positive constant $M$ such that 
\begin{equation}
\nonumber 
\left|\frac{a_m\lm + b_m}{c_m\lm+d_m}\right|\leq M, \ \ \
\left|\frac{a_m}{c_m}\right|\leq M.
\end{equation}
Hence we have also
\begin{equation}\label{approx3}
\lim_{m\to\infty}|a_m\lm+b_m|=0
\end{equation}
and 
\begin{equation}\label{approx4}
\lim_{m\to\infty}|a_m(c_m\lm+b_m)^2|=0.
\end{equation}
Put
\[
\ttau_*(g_m) = \left(
\begin{array}{cc}A_m & B_m \\ C_m & D_m\end{array}
\right).
\]

\begin{lem}\label{CLD} \ 
$\lim_{m\to\infty}\|C_mL_\lm + D_m\|=0$.
\end{lem}
{\bf Proof} \  
We calculate the components of $C_mL_\lm+D_m$. 
Put
$C_mL_\lm + D_m 
= \MAT{\al_{11}}{\al_{12}}{\al_{21}}{\al_{22}}. 
$
Then we have
\begin{eqnarray*}
\label{a11}
\al_{11} &=& a_m(c_m\lm+d_m)^2+2(a_m\lm+b_m)c_m(c_m\lm+d_m) \\
\label{a12}
\al_{12} &=& (a_m\lm+b_m)(c_m\lm+d_m)^2
- a_m\lm(c_m\lm+d_m)^2 
- 2(a_m\lm+b_m)c_m\lm(c_m\lm+d_m) \\
\label{a21}
\al_{21} &=& 3c_m(c_m\lm+d_m)^2 \\
\label{a22}
\al_{22} &=& (c_m\lm+d_m)^3 - 3c_m\lm(c_m\lm+d_m)^2.
\end{eqnarray*}
Then 
\[ \lim_{m\to\infty}\al_{ij}=0 \]
follows easily from $(\ref{approx2}), (\ref{approx3})$, and $(\ref{approx4})$. 
\ \qed 

\bigskip
{\bf Proof of the proposition (continued)}. \ 
Suppose that $\ell_\lm \cap \bDelta$
contains a point $a=[a':a''] \in \P^3$, where 
$a'=L_\lm a''$. Then, we have
\[ \|(C_gL_\lm+ D_g)a''\| \geq \|a''\| \]
for any $g \in \Gamma$.  Since $a''\neq 0$,
this contradicts Lemma \ref{CLD}. 
\ \qed 

\begin{rem} The condition that $\Gamma$ should not contain 
parabolic elements is indispensable. Indeed, the group
$\tGamma$ induced by the rank 2 abelian group
$\Gamma=\{\tau_1, \tau_2\}$, $\tau_1(z)=z+1$, $\tau_2(z)=z+i$, 
gives a counter example.
\end{rem}

\subsubsection{Kleinian groups acting on a quadric surface}
Let $S \in \P^3$ be the quartic surface $S : z_0z_3 - z_1z_2 = 0$, and 
\[  q : \P^1 \times \P^1 \to S \]
the Segre map $q([u_0:u_1], [v_0:v_1]) = [u_0v_0:u_0v_1:u_1v_0:u_1v_1]$. 
We consider the case where a subgroup $\Gamma \subset \PSL_2(\C)$
acts trivially on the second component of $\P^1\times \P^1$. 
This case was studied by Fujiki \cite{Hayama} and Guillot \cite{Gui}. 
Here, we shall reprove a theorem of Fujiki, as an application of 
Theorem \ref{Ford_Rgn}.

Then the Segre map $q$ defines a group representation
\[ q_* : \Gamma \to \PGL_4(\C), \]
which is induced by the following commutative diagram:
\[
\begin{array}{rcl}
                \P^1 \times \P^1 &  \stackrel{q}{\longrightarrow} & \P^3  \\
g\times 1 \downarrow \hspace{1.5em}&      & \downarrow  q_*(g) \hspace{1em} \\
                \P^1 \times \P^1 &  \stackrel{q}{\longrightarrow} & \P^3.
\end{array}
\]
Explicitly, for 
$g=\pm\left(
\begin{array}{cc}
a & b \\ c & d
\end{array}
\right) \in \PSL_2(\C)$,  
$\tg=q_*(g)$ is given by 
\begin{equation}\label{tgF}
\tg = \pm\left(
\begin{array}{cc}
aI  & bI \\
cI & dI 
\end{array}
\right) \in \PGL_4(\C),   
\end{equation}
where $I$ denote the identity matrix of size $2$.

\bigskip
\underline{\bf Limit sets}  \ 
Let $\Gamma \subset \PSL_2(\C)$ be a Kleinian group whose 
set of discontinuity $\Omega_{\P^1}$ contains $[1:0]\in\P^1$.
Put $\Lambda_{\P^1}=\P^1\setminus\Omega_{\P^1}$
and $\tGamma = q_*(\Gamma)$. 
\begin{prop}\label{one_side_lim}
The limit set of $\tGamma$ is given by 
\[
\Lambda(\tGamma) = q(\Lambda_{\P^1}\times\P^1). 
\]
Thus $\tGamma$ is of type $\ELL$, 
and satisfies $(\clubsuit)$ and $(\spadesuit)$. 
\end{prop}
{\bf Proof} \  
As in subsection \ref{tw}, we embed $\Grass$ into 
$\P^5=\P(\wg^2\C^4)$, and 
consider the group homomorphism
\[ \bar\rho : \PGL_4(\C) \to \PGL_6(\C). \]
Let $g=\pm\MAT{a}{b}{c}{d}\in\PSL_2(\C)$. 
With respect to the basis 
\[ 
\{ e_0\wg e_1, \ e_0\wg e_2, \ e_0\wg e_3, 
  \ e_1\wg e_2, \ e_1\wg e_3, \ e_2\wg e_3\},
\]
of $\wg^2\C^4=\C^6$, the matrix 
\begin{equation}\label{size6d}\nonumber
G(a,b,c,d):=\left(
\begin{array}{rrrrrr}
a^2 & 0 & ab &  -ab & 0 & b^2 \\
0 & 1 & 0 & 0 & 0   & 0 \\
ac & 0 & ad & -bc & 0 & bd \\
-ac & 0 & -bc & ad & 0 & -bd \\
0 & 0 & 0 & 0 & 1 & 0 \\
c^2 & 0 & cd & -cd & 0 &d^2
\end{array}
\right) \in \SL_4(\C)
\end{equation}
represents $\bar\rho(\tg) \in \PGL_6(\C)$.
Let $\{\gamma_n\}\subset\Gamma$ be a normal sequence.
Let 
\[ 
g_n = \MAT{a_n}{b_n}{c_n}{d_n} \in \SL_2(\C), \ \ \ n=1,2,\dots 
\]
be a sequence of representatives of $\{\gamma_n\}$.
Since $[1:0]\in\Omega_{\P^1}$, 
$\{c_n^{-1}g_n\}$ converges to a matrix of the form 
$h=\MAT{\mu}{-\lm\mu}{1}{-\lm}\in M_2(\C)$, $\lm,\mu\in\C$. 
Letting $G_n=c_n^{-2}G(a_n,b_n,c_n,d_n)$, 
we calculate the limit : 
\begin{small}
\begin{equation}\label{rhogamma_d}\nonumber
G:=\lim_n G_n=\left(
\begin{array}{rrrrrr}
\mu^2 & 0 & -\lm\mu^2 &  \lm\mu^2 & 0 & \lm^2\mu^2 \\
0 & 0 & 0 & 0 & 0 & 0 \\
\mu & 0 & -\lm\mu & \lm\mu & 0 & \lm^2\mu \\
-\mu & 0 & \lm\mu & -\lm\mu & 0 & -\lm^2\mu \\
0 & 0 & 0 & 0 & 0 & 0 \\
1& 0 & -\lm & \lm & 0 & \lm^2
\end{array}
\right). 
\end{equation}
\end{small}
Thus $G$ defines a projection to a single point,  
\[ \P^5\setminus H \to I=\{[\mu^2:0:\mu:-\mu:0:1]\}, \]
where $H$ is the 4-plane defined by 
\[ 
H = \{ \zeta \in \P^5 : 
\zeta_0-\lm\zeta_2 +\lm\zeta_3 + \lm^2\zeta_5=0 \}. 
\]
Note that $I$ is contained in $\Grass$ and corresponding 
to the line 
\begin{equation}\label{LL}
 (\mu e_0+ e_2)\wg(\mu e_1 + e_3)=\{z'=\mu z''\} 
\end{equation}
in $\P^3$. This line coincides with $q([\mu:1]\times \P^1)$.
That $\tGamma$ satisfies $(\clubsuit)$ and $(\spadesuit)$ follows 
from the form $(\ref{tgF})$ and the fact 
that $\sum_{g\in\Gamma^*}|c_g|^{-4}<+\infty$ 
in the Kleinian group theory (\cite[Theorem II.B.5]{Maskit}). 
\ \qed

\bigskip
\underline{\bf Compact Quotients} \ 
As an application of Theorem \ref{Ford_Rgn}, we have the following. 
\begin{thm}\label{Hayama}\cite{Hayama} 
If $\Gamma$ is convex-cocompact, 
then $\left(\P^3\setminus \Lambda(\tGamma)\right)/\tGamma$ 
is compact. 
\end{thm}
{\bf Proof} \  
The outline of the proof is the same as that of Theorem \ref{TQ}. 
As in that proof, it is sufficient to prove the following. 
\ \qed

\begin{prop} Any limit line does not intersect $\bar{\Delta}$.
\end{prop}
{\bf Proof} \   A limit line $\ell_\lm$ is given by 
$z'=\lm z''$ by (\ref{LL}), where $[\lm:1]\in\P^1$ is 
the limit point of $\Gamma$. 
Now, suppose that there exits a limit line $\ell_\lm$ such that
$\ell_\lm\cap\bar\Delta$ is non-empty. Take a point
 $[a':a'']\in\ell_\lm\cap\bar\Delta$. Then, by $a'=\lm a''$,
$\tg=\MAT{a_gI}{b_gI}{c_gI}{d_gI}$, and 
$\|C_g a'+D_g a''\|\geq\|a''\|$,
we have $\|(c_g\lm+d_g)a''\|\geq\|a''\|$ for any $g\in\Gamma$.
Since $a''\neq 0$, this contradicts (\ref{approx2}). 
\ \qed

\subsubsection{Kleinian groups acting on a cone over a conic}
For the moment we have only a very simple example of type $\ELL$ in this case.
Many discrete subgroups acting on the cone can be constructed by the method used 
in \cite[p.278]{plane}. It is plausible some of them are of type $\ELL$, but their 
canonical quotients will be non-compact.

\noindent
{\bf Example} \ Let $\Gamma\subset \SL_4(\C)$ be an infinite cyclic group generated by
\[
g=\left(\begin{array}{cccc}
\al^2 & 0 & 0 & 0 \\
0 & 1 & 0 & 0 \\
0 & 0 & \al^{-2} & 0 \\
p & q & r & 1 \end{array}\right), \ \ |\al|>1.
\]
With respect to the basis 
\[
\{e_0\wg e_1,  \ e_0\wg e_2,  \ e_1\wg e_2, \ e_0\wg e_3, \ e_1\wg e_3, \ 
e_2\wg e_3 \},
\]
we have 
\[
\rho(g^n)=\left(\begin{array}{cccccc}
\al^{2n} & 0 & 0 & 0 & 0 & 0 \\  
0 & 1 & 0 & 0 & 0 & 0 \\  
0 & 0 & \al^{-2n} & 0 & 0 & 0 \\  
n\al^{2n}q & \frac{1-\al^{2n}}{\al^{-2}-1}r & 0 & \al^{2n} & 0 & 0 \\  
-\frac{\al^{2n}-1}{\al^2-1}p  & 0 &\frac{\al^{-2n}-1}{\al^{-2}-1}r & 0 & 1 & 0 \\  
0 & -\frac{1-\al^{-2n}}{\al^2-1}p & -n\al^{-2n}q & 0 & 0 & \al^{-2n}
\end{array}\right)
\in \PGL_6(\C).
\]
This implies that the limit image of the sequence $\{\rho(g^n)\}$, \ $n \to +\infty/-\infty$,
is a point if and only if $q\neq 0$. If $q\neq 0$, there are exactly two limit lines, which are
\[
\ell_1 = e_0\wg e_3 \ \ \mbox{and} \ \ \ell_2=e_2\wg e_3 . 
\]
Thus $\Gamma$ is of type $\ELL$ if and only if $q\neq 0$. The cone
$S=\{z_0z_2-z_1^2=0\}$ contains $\ell_1$ and $\ell_2$, and
they are invariant by $\Gamma$. Note that the quotient space 
$\Omega(\Gamma)/\Gamma=(\P^3\setminus\{\ell_1\cup\ell_2\})/\Gamma$
contains a non-compact surface 
$(S\setminus \{\ell_1\cup\ell_2\})/\Gamma$ as a closed submanifold, 
which is a $\C$-bundle over the elliptic curve $\C^*/\langle\al\rangle$. 
Therefore $\Omega(\Gamma)/\Gamma$ is not compact. The group satisfies
$(\clubsuit)$ and $(\spadesuit)$ with respect to a new system $[w]$ 
of coordinates, such as 
$w_0=z_0+z_1-z_2-z_3, \ w_1=-z_0+z_1+z_2-z_3, \ w_2=z_0-z_1+z_2+z_3, \ w_3=z_3$.

\begin{flushleft}
Sophia University\\
Kioicho 7-1, Chiyoda-ku, Tokyo 102-8554, \ JAPAN \\
email : masahide.kato@sophia.ac.jp
\end{flushleft}
\end{document}